\documentclass[11pt,reqno]{amsart}
\usepackage{amsmath,amssymb,amsthm,url}
\newtheorem{theorem}{Theorem}
\newtheorem{corollary}{Corollary}
\newtheorem{lemma}{Lemma}

\theoremstyle{definition}
\newtheorem{definition}{Definition}
\theoremstyle{remark}
\newtheorem{remark}{Remark}
\newtheorem{example}{Example}
\title[V-bounded solutions]
{Existence of V-bounded solutions\\ for nonautonomous nonlinear systems\\ via
the Wa\.zewski topological principle}
\author[V. Lagoda]
{Volodymyr Lagoda}
\author[I. Parasyuk]
{Igor Parasyuk} \email{pio@univ.kiev.ua}
\address{National Taras Shevchenko Univesity of
Kyiv\\ Volodymyrs'ka 64, Kyiv, 01033, Ukraine}

\begin{document}

\begin{abstract}{We establish a number of new sufficient conditions for the existence of global
(defined on the entire time axis) solutions of nonlinear nonautonomous systems
by means of the Wa\.zewski topological principle. The systems under
consideration are characterized by the monotonicity property with respect to a
certain auxiliary  guiding function $W(t,x)$ depending on time and phase
coordinates. Another auxiliary function $V(t,x)$, such that $\lim_{\|x\|\to
\infty}V(t,x)= \infty$ for all $t\in \mathbb R $, is used to estimate the location of global
solutions in the extended phase space. The approach developed is applied to
Lagrangian systems, and in particular, to establish new sufficient conditions
for the existence of almost periodic solutions.}
\end{abstract}
\maketitle

\section{Introduction}
This paper is a modified and extended version of our e-print \cite{LagPar09}.
 Its goal is to lay down some new sufficient conditions
 under which the nonlinear nonautonomous system of ODEs
\begin{equation}\label{eq:nlsys} \dot{x}=f(t,x),
\end{equation}
where $f:\Omega  \mapsto \mathbb R^n$ ($\Omega\subseteq \mathbb R^{1+n}$),  has a global solution
$x(t)$ which exists on the entire time axis and  possess the property that a
given auxiliary spatially coercive  function $V(t,x)$ (a time dependent norm
surrogate) is bounded along the graph of $x(t)$. We especially focus on getting
estimates for the function $V(t,x(t))$. The main results are obtained by using
the Wa\.zewski topological principle \cite{Waz47,Har70,Con76,Con78}, and some
of them generalize  results of V.~M.~Cheresiz \cite{Che74}.

It should be noted that the Wa\.zewski topological principle was successfully
exploited for proving the existence of bounded solutions to some boundary value
problems in \cite{Con75} and to  quasihomogeneous systems in
\cite{Iva85a,Iva85b} (see also a discussion in \cite{Ort08}).

In order to apply the the Wa\.zewski principle, along with the function $V(\cdot)$
we use another auxiliary function $W(t,x)$ with positive derivative by virtue
of the system (\ref{eq:nlsys}) in the domain where $V>0$. We call $V$ and $W$
\emph{the estimating function} and \emph{the guiding function} respectively and
we say that together they form the V--W-pair of the system. Note that the term
''guiding function'' we borrow from \cite{KraZab75} (originally --- ''guiding
potential''). Basically topological method of guiding functions, which was
developed by M.~A.~Krasnosel'ski and A.~I.~Perov, is an effective tool for
proving the existence of bounded solutions of essentially nonlinear systems too
(see the bibliography in \cite{KraZab75, MawWar02}). But, except
\cite{Ort08,Avr03}, in all papers known to us, only independent of time guiding
functions were used.

In \cite{Che74}, the role of V--W-pair plays some function of Euclidean norm
together with an indefinite nondegenerate quadratic form. It appears that in
this case sufficient conditions for the existence of bounded solutions as well
as the estimates of their norms coincide with those obtained by means of
technique developed in \cite{PerTru74, PerTru86} for indefinitely monotone (not
necessarily finite dimensional) systems.

We shall not mention here another interesting approaches in studying the
existence problem of bounded solutions to nonlinear systems, because they have
not been used in this paper. For the corresponding information the reader is
referred to \cite{MNS72,Per03,
War94,SBB02,Blo88,BerZha95,BCM97,ZakPar99,Cie03,Cie03b,CheMam05,Her06,
ARS06,Sly02, Gro06,Kar07}.

This paper is organized as follows. Section~\ref{defVWpair} contains necessary
definitions, in particular, the notion of V--W-pair is introduced and some
additional conditions imposed on estimating and guiding functions are
described. In section~\ref{ExistUniqVBS}, we prove two main theorems concerning
the existence and the uniqueness of V-bounded solution to a nonlinear
nonautonomous system possessing V--W-pair.  In section~\ref{QuadForms} we show
how the results of section~\ref{ExistUniqVBS} can be applied in the case where
the estimating and guiding functions are constructed by means of nonautonomous
quadratic forms. In this connection it should be pointed out that guiding
quadratic forms play an important role in the theory of linear dichotomous
systems with (integrally) bounded coefficients \cite{DalKre70,MSK03,Sam02}. As
an example of application of our technique, we generalize results of
\cite{MNS72,Per03} on the existence of bounded solutions to quasilinear
nonautonomous system with exponentially dichotomic linear part. Finally, in
section~\ref{Lagrsys}, the approach developed in section~\ref{ExistUniqVBS} is
applied to a quasiconvex Lagrangian system of mechanical type with time-varying
holonomic constraint. For such systems, we establish sufficient conditions for
the existence of global solutions along which the Lagrangian function remains
bounded. The case of almost periodic Lagrangian is also discussed. As an
example we consider motion of a particle on helicoid under the impact of force
of gravity and repelling potential field of force. Note that bounded and almost
periodic solutions of globally strongly convex and Lipschitzian Lagrangian
systems were studied in~\cite{Cie03b}.

\section{The definition of V--W-pair
 and the main assumptions} \label{defVWpair}

Let  $\Omega$ be a domain in $\mathbb R^{1+n}=\{t\in \mathbb R\}\times \{x\in \mathbb R^{n}\}$ such that
the projection of $\Omega $ on the time axis $\{t\in \mathbb R\}$ covers all this axis,
and let  $f(\cdot)\in \mathrm{C}(\Omega \mapsto \mathbb R^n)$. It will be always assumed that each
solution of the system \eqref{eq:nlsys} has the uniqueness property.

\begin{definition}A function $V(\cdot)\in  \mathrm{C}(\mathbb R\times \mathbb R^n\! \mapsto\! \mathbb R)$
of variables $t\in \mathbb R,\;x\in \mathbb R^{n}$ will be called spatially coercive, if for
any $t\in \mathbb R$ the function $V_t(\cdot):=V(t,\cdot):\mathbb R^n \mapsto \mathbb R$ has the following
properties: the level set $V_t^{-1}(0):=\{x\in \mathbb R^n:V_t(x)=0\}$ is nonempty and
$\lim_{\|x\|\to \infty}V_t(x)= \infty$. If in addition $V(\cdot)\in  \mathrm{C}^1(\mathbb R\times \mathbb R^n\!
\mapsto\! \mathbb R)$ and $\|\frac{\partial V_t(x)}{\partial x}\|> 0 $ once $V_t(x)> 0$, then $V(\cdot)$
will be called a regular spatially coercive function.
\end{definition}

Note that for each $t\in \mathbb R$ and each $v>0$ the level set $V_t^{-1}(v)$ of
regular spatially coercive function $V(\cdot)$ is a compact connected and simply
connected hypersurface which, thus, is homeomorphic to $(n-1)$-dimensional
sphere;  in addition, if $v_2>v_1$, then the set
$V_t^{-1}\left((-\infty,v_1]\right):=\{x\in \mathbb R^{n}: V_{t}(x)\le v_1\}$ is a proper subset
of the set $V_t^{-1}\left((-\infty,v_2]\right)$.

\begin{definition}For a spatially coercive function $V(\cdot)$, a global solution $x(t)$, $t\in I$
of the system (\ref{eq:nlsys}) is said to be V-bounded if  $$\sup_{t\in
I}V(t,x(t))<\infty, $$ and $V(\cdot)$ is then called an estimating function.
\end{definition}

For any $U(\cdot)\in \mathrm{C}^1(\Omega  \mapsto \mathbb R)$, define  $$\dot{U}_{f}:=\frac{\partial U}{\partial
t}+\frac{\partial U}{\partial x}\cdot f.$$

\begin{definition} A function $W(\cdot)\in \mathrm{C}^1(\Omega \! \mapsto\! \mathbb R)$
will be called a guiding function concordant with a spatially coercive function
$V(\cdot)$ if $\Omega \cap V^{-1}\left((0,\infty)\right)\ne \varnothing$ and $
 \dot W_f(t,x)>0$ for all $(t,x)\in  \Omega \cap V^{-1}\left((0,\infty)\right)$.
  \end{definition}

\begin{definition} A regular spatially coercive function $V(\cdot)$
together with a concordant  guiding function $W(\cdot)$
 will be called a V--W-pair of the system \eqref{eq:nlsys}.
\end{definition}

Denote by $\Pi_t:=\{t\}\times \mathbb R^{n}$ the ''vertical'' hyperplane in
$\mathbb R^{1+n}=\mathbb R\times \mathbb R^{n}$ passing through $(t,0)$ and for any set $\mathcal{A}\subset
\mathbb R\times  \mathbb R^{n}$ denote by $\mathcal{A}_t$ the natural projection of the set $\Pi_t\cap
\mathcal{A}$ onto $\mathbb R^n$.

In so far, we suppose that the system (\ref{eq:nlsys}) has a V--W-pair which
satisfies the following additional conditions:
\medskip
\begin{itemize}
\item[(A):] there exist numbers
 $w^*,\;w_*$ ($w^*>w_*$), $c^*>0, c_*\in [0,\infty]$, and a connected component $\mathcal{W}$ of the set
 $W^{-1}\left((w_*,w^*)\right)$ such that
 for any $t\in \mathbb R$ the number $w^*$ belongs to the range of function
$W_{t}(\cdot):=W(t,\cdot):\Omega_t \mapsto \mathbb R$, the set $V^{-1}\left((-\infty,0]\right)$ belongs to
$\mathcal{W}$, and the following inequalities hold
\begin{equation}\label{eq:mainineq}
-c_* \dot W_f(t,x)\le \dot V_f(t,x)\le c^*\dot W_f(t,x)\quad \forall(t,x)\in
  V^{-1}\left([0,\infty)\right)\cap \mathcal{W};
\end{equation}
\end{itemize}

Note that from condition~(A), it follows that
\begin{equation}\label{eq:w0w0}
\begin{split}
  w_0(t) & :=\min\{W_t(x):x\in V_t^{-1}(0)\} > w_*, \\
   w^0(t) & :=\max\{W_t(x): x\in V_t^{-1}(0)\}
  <  w^*,
\end{split}
\end{equation}
thus, the set $\partial \mathcal{W}\cap W^{-1}(w^*)$ coincides with the set of  exit points
from $\mathcal{W}$, each point of $\partial \mathcal{W}\cap W^{-1}(w^*)$ being a strict exit point.
Denote
\begin{gather*}
  \mathcal{W}^{se}:=\partial \mathcal{W}\cap W^{-1}(w^*).
\end{gather*}

\begin{itemize}
  \item[(B):] the function
\begin{gather*}
  \alpha(t):=\inf \left\{\dot W_f(t,x):x\in V_t^{-1}\left((0,\infty)\right)\cap \mathcal{W}_t\right\}
\end{gather*}
has the property
  $\int_{-\infty}^{0}\alpha (s)\,\mathrm{d}s=\int_{0}^{\infty}\alpha (s)\,\mathrm{d}s=\infty$;

\item[(C):]  for any sufficiently large by absolute value negative  $t$, the
Wa\.zewski condition is fulfilled: there exists a bounded subset $\mathcal{M}_t$ of
the set $\mathcal{W}_t\cup\mathcal{W}_t^{se}$  such that the set $ \{t\}\times \left(\mathcal{M}_{t}\cap
\mathcal{W}_t^{se}\right)$ is a retract of $\{(s,x)\in \mathcal{W}^{se}:s\ge t\}$, but is not a
retract of $\{t\}\times \mathcal{M}_t$.

\end{itemize}

\begin{remark}\label{rem:1}
In the case where $V(\cdot)$ and $W(\cdot)$ do not depend of $t$, one can consider the
inequalities \eqref{eq:mainineq} as an analogue of  the regularity condition
for the guiding function  $W(\cdot)$ (see \cite{KraZab75}). The main consequence of
regularity in this case is that the pair $cW(\cdot)$, $c W(\cdot)-V(\cdot)$ (or $cW(\cdot)$,
$cW(\cdot)+V$) is a complete set of guiding functions for any $c>c^*$ (for any
$c>c_*$ if $c_*<\infty$).
\end{remark}

\begin{remark} The condition~(C) is fulfilled if for any negative
sufficiently large by abso\-lute value  $t$ there exists a  compact manifold
$\mathcal{M}_t$ with border $\partial \mathcal{M}_t$ such that the interior of $\mathcal{M}_{t}$ belongs
to $\mathcal{W}_t$ and the set $ \{t\}\times \left(\mathcal{M}_{t}\cap \mathcal{W}_t^{se}\right)$ is a retract
of $\{(s,x)\in \mathcal{W}^{se}:s\ge t\}$.  In fact, as is well known,  $\partial \mathcal{M}_t$ is
not a retract of $\mathcal{M}_t$.

Taking into account that $\mathcal{W}^{se}$ is  a  connected component of regular
level hypersurface $W^{-1}(w^*)$, the condition~(C) can be replaced by the
following weaker condition:
\end{remark}

\begin{itemize}
  \item[(C$'$):]  there exists a bounded subset $\mathcal{M}_t$ of
the set $\mathcal{W}_t\cup\mathcal{W}_t^{se}$ which cannot be continuously imbedded into
$\mathcal{W}^{se}$ in such a way that the image of $\mathcal{M}_{t}\cap \mathcal{W}_t^{se}$ is $
\{t\}\times \left(\mathcal{M}_{t}\cap \mathcal{W}_t^{se}\right)$.
\end{itemize}

\section{The existence and the uniqueness
 of V-bounded solution}\label{ExistUniqVBS}

The lemma given below open the door to estimation of solutions of the system
(\ref{eq:nlsys}) by means of V--W-pair.

\begin{lemma}\label{lem:estimate}
Suppose that the system (\ref{eq:nlsys}) has V--W-pair satisfying the
condition~(A). Let  $x(t)$ be  such a solution of (\ref{eq:nlsys}) that
$(t,x(t))\in \mathcal{W}$ for all $t\in [t_0,t_1]$.

Then the following assertions are true:

-- if   $V(t,x(t))>0$ for all $t\in (t_0,t_1]$, then
\begin{gather}\label{eq:estV1}
  V(t,x(t))\le V(t_0,x(t_0))+ c^*\left[w^* - W(t_0,x(t_0))\right]\quad \forall t\in
  [t_0,t_1];
\end{gather}

-- if $V(t,x(t))>0$ for all $t\in (t_0,t_1)$ and $V(t_1,x(t_1))=0$, then
\begin{gather}\label{eq:estV2}
  V(t,x(t))\le
    \frac{c_*}{c_*+c^*}V(t_0,x(t_0))+\frac{c_*c^*}{c_*+c^*}\left[w^0(t_1)-W(t_0,x(t_0))\right]
     \\ \forall t\in [t_0,t_1];\nonumber\end{gather}

-- if the condition~(B) is fulfilled, $V(t_0,x(t_0))\ge0$  and
\begin{gather} \label{eq:esttau}
  \int_{t_0}^{t_1}\alpha (s)\,\mathrm{d}s\ge w^*-w_*,
\end{gather}
then there exists  $\tau  \in  (t_0,t_1)$ such that $V(\tau,x(\tau))= 0$.
\end{lemma}
\begin{proof} Let the condition~(A) is fulfilled. Put $v(t):=V(t,x(t))$,
$w(t)=W(t,x(t))$. The inequality \eqref{eq:estV1} obviously follows from $\dot
v(t)\le c^*\dot w(t)$, $t\in [t_0,t_1]$. In order to prove the inequality
\eqref{eq:estV2}, denote by $\hat t$ any point where $v(t)$ reaches its maximum
on $[t_0,t_1]$ and observe that
\begin{gather*}
  w(t_1)-w(t_0)=\int_{t_0}^{\hat t}\dot w(t)\,\mathrm{d}t +\int_{\hat t}^{t_1}\dot w(t)\,\mathrm{d}t\ge
\frac{1}{c^*}\int_{t_0}^{\hat t}\dot v(t)\,\mathrm{d}t-\frac{1}{c_*}\int_{\hat t}^{t_1}\dot
v(t)\,\mathrm{d}t=\\ \frac{v(\hat t)-v(t_0)}{c^*}+\frac{v(\hat
t)}{c_*}\ge\frac{(c_*+c^*)v(t)}{c_*c^*}-\frac{v(t_0)}{c_*}.
\end{gather*}
Since $v(t_1)=0$, then $w(t_1)\le w^0(t_1)$ and we get \eqref{eq:estV2}.

Now let the condition~(B) is fulfilled and $v(t_0)\ge 0$. If we assume that
$v(t)>0$ on $(t_0,t_1)$, then
\begin{gather*}
  \int_{t_0}^{t_1}\alpha(t)\,\mathrm{d}t\le \int_{t_0}^{t_1}\dot
  w(t)\,\mathrm{d}t=w(t_1)-w(t_0)<w^*-w_*.
\end{gather*}
This contradicts the inequality \eqref{eq:esttau}.
\end{proof}

Put
\begin{gather}
   \omega_0:=\inf_{t\in \mathbb R}w_0(t),\quad \omega^0:=\sup_{t\in \mathbb R}w^0(t),\label{eq:defom0om0}\\
  \nu :=\liminf_{t\to - \infty}\,\sup \left\{V_t(x)-c^*W_t(x):x\in
  \mathcal{M}_t\cap V_t^{-1}\left((0,\infty)\right)\right\}.\label{eq:defnu}
\end{gather}

Now we are in position to prove the following theorem.

\begin{theorem}\label{th:exbs} Assume that the system (\ref{eq:nlsys})
has a V--W-pair satisfying the conditions~(A),(B),(C) (or (C\,$'$)),  and $\nu
<\infty $. Let there exists a number $V^*>c^*w^0+\max\left\{\nu,-c^*\omega_0\right\}$ such
that $$\mathrm{cls}\left(V^{-1}\left([0,V^*)\right)\cap\mathcal{W}\right)\subset \Omega $$ (here $\mathrm{cls}$
means the closure operation). Then the system (\ref{eq:nlsys}) has a V-bounded
global solution $x_\ast  (t),\;t\in \mathbb R$, which for all $t\in \mathbb R$
 satisfies the inequalities
\begin{gather}\label{eq:estVfin}
  V(t,x_*(t))\le \frac{c_*c^*}{c_*+c^*}\left[\sup_{t\le s\le \tau_+(t)}w^0(s)-
  \inf_{\tau_-(t)\le s\le t}w_0(s)\right]
  \le \\ \frac{c_*c^*}{c_*+c^*}(\omega^0-\omega_0),\nonumber
     \\  \omega_0\le W(t,x_*(t))\le \omega^0  \end{gather}
where $\tau_+(t)$ and $\tau_-(t)$ are, respectively, the roots of equations
\begin{gather*}
  \int_{t}^{\tau_+}\alpha(s)\,\mathrm{d}s=\omega^0-\omega_0,\quad \int_{\tau_-}^{t}\alpha(s)\,\mathrm{d}s=\omega^0-\omega_0.
\end{gather*}
\end{theorem}

\begin{proof}
Firstly observe that we may consider the numbers $w^*\in (\omega^0,\infty)$ and
$w_*\in (-\infty,\omega_0)$ to be arbitrarily close to $\omega^0$ and $\omega_0$ respectively.
 From definitions of $\nu
$ and $V^*$ it follows that there exists a sequence $t_j\to -\infty,\;j\to \infty $,
such that
\begin{gather}\label{eq:tjV*}
 c^*w^*+ \sup \left\{V_{t_j}(x)-c^*W_{t_j}(x):x\in
  \mathcal{M}_{t_j}\cap V_{t_j}^{-1}\left((0,\infty)\right)\right\}< V^*.
\end{gather}
In view of condition~(C) (or (C$'$)) from Wa\.zewski principle it follows that
for any $j$ there exists a point $x_{0j}\in \mathcal{M}_{t_j}$ such that the global
solution $x_j(t),\;t\in I_j$, which satisfies the initial condition
$x_j(t_j)=x_{0j}$ has the property
\begin{gather*}
  (t,x_j(t))\in \mathcal{W}  \quad
  \forall t\in [t_j,\infty)\cap I_j.
\end{gather*}
Let us show that
\begin{equation}\label{eq:estvjt}
  v_j(t):=V(t,x_j(t))<V^*\quad  \forall t\in I_j\cap[t_j,\infty).
\end{equation}

For any natural $j$, it is sufficient to consider the following cases: (I)
$v_j(t)> 0$ for all $t\in I_j\cap(t_j,\infty)$; (II) $v_j(t_j)\ge 0$, there exist
$t_*\ge t_j$ and $t^*\ge t_*$ such that $v(t_*)=v(t^*)=0$, $v(t)> 0$ for all
$t\in I_j\cap(t^*,\infty)$, and if $t_*>t_j$, then $v_j(t)>0$ for all $t\in
(t_j,t_*)$; (III) there exist increasing sequences $t_{k\ast},t^{\ast}_{k}$ in
$I_j\cap[t_j,\infty)$, $k\in \mathbb N$,  such that $t_{k\ast}< t^*_{k}$, $t_{k+1,*}\ge
t^{\ast}_k$, $t^{\ast}_{k}\to \sup\left\{t\in I_j\right\}$, $k\to \infty,$ and
\begin{gather*}
  v_j(t_{k\ast})=v_j(t^{\ast}_{k})=0,\\
   v_j(t)>  0\;\forall t\in (t_{k\ast},t^*_{k}), \quad  v_j(t)\le 0
\;\forall t\in (I_j\setminus \textstyle\bigcup_{k=1}^\infty(t_{k\ast},t^{\ast
}_k)\cap[t_j,\infty).
\end{gather*}

In the case (I), observe that for sufficiently small $\delta >0$ we have
$$v_j(t_j)+c^*\left[w^*-W(t_j,x_{0j})\right] <c^*w^*+\nu +\delta <V^*.$$ Now the inequality
\eqref{eq:estvjt} immediately  follows from \eqref{eq:estV1}.

In the case (II), observe that
\begin{gather*}
  v(t^*)+c^*\left[w^*-W(t^*,x_j(t^*))\right]\le
  c^*\left[w^*-w_0(t^*)\right]\le c^*\left[w^*- \omega_0\right] <V^*.
\end{gather*}
Thus, similarly to the case (I), we obtain the estimate $v_j(t)< V^*$ for all
$t\in I_j\cap[t^*,\infty)$. Next, if $t_*>t_j$, then $W(t_j,x_j(t_j))\le
W(t_*,x_j(t_*))\le w^0(t_*)$ and from \eqref{eq:estV2} it follows that
\begin{gather*}
  v_j(t)\le
  \frac{c_*}{c_*+c^*}v_j(t_j)+\frac{c_*c^*}{c_*+c^*}\left[w^0(t_*)-W(t_j,x_j(t_j))\right]<
  \frac{c_*}{c_*+c^*}V^*\le V^* \\
  \forall t\in [t_j,t_*].
\end{gather*}
If now $t_*=t^*$, then the inequality \eqref{eq:estvjt} holds true. And if
$t_*<t^*$, then for any successive zeroes $t_1,t_2\in [t_*,t^*]$ of function
$v_j(t)$  from \eqref{eq:estV2} it follows that
\begin{gather*}
  v_j(t)\le \frac{c_*c^*}{c_*+c^*}\left[w^0(t_2)-w_0(t_1)\right]<\frac{c_*c^*}{c_*+c^*}\left[\omega^0-\omega_0\right]
  <V^*\quad \forall t\in
  [t_1,t_2].
\end{gather*}
Thus, we obtain inequality \eqref{eq:estvjt} in the case (II), and now  it
becomes obvious that this inequality is valid also for the case (III).

The above reasoning allows us to make conclusion that in view of definition of
$V^*$ the graph of $x_j(t)$, $t\in I_j\cap [t_j,\infty)$, is contained in a closed
subset of $\mathcal{W}$. This yields inclusion $[t_j,\infty)\subset I_j$.

Now we are in position to prove the existence of V-bounded solution $x_*(t)$ by
the known scheme (see, e.g., \cite{Che74,Iva85b,KraZab75}). Namely, if we
denote by $x(t,t_0,x_0)$  the solution  which for $t=t_0$ takes the value
$x_0$, then setting $\xi_j:=x_j(0)$, we obtain the equalities
\begin{gather*}
  x_j(t)=x(t,0,x_j(0))=x(t,0,\xi_j),\quad t\in [t_j,\infty).
\end{gather*}
Having selected from the sequence $\xi_j \in
\mathrm{cls}\,\left(V_0^{-1}([0,V^*])\cap\mathcal{W}_0\right)\subset \Omega_0 $ a subsequence
converging to $x_*$, put $x_*(t):=x(t,0,x_*)$.  Using reductio ad absurdum
reasoning it is easy to show that on the maximal existence interval $I$ of this
solution we have the inclusion $$(t,x_*(t))\in \mathrm{cls}(V^{-1}\left([0,V_*)\right)\cap
\mathcal{W}).$$ Therefore  $I=\mathbb R$.

Now we are able to establish a sharper estimate for $v_*(t):=V(t,x_*(t))$.
Namely, for any $t\in \mathbb R$ such that $v_*(t)>0$, in virtue of
Lemma~\ref{lem:estimate}, the point $t$ lies between two successive zeroes
$t_*(t),t^*(t)$ of $v_*(t)$ each of which is contained in the segment
$[\tau_-(t),\tau_+(t)]$. Then the inequality \eqref{eq:estVfin} easily follows from
\eqref{eq:estV2} once we put there $t_0=t_*(t),\;t_1=t^*(t)$.
\end{proof}

The following theorem establishes sufficient conditions for the uniqueness of
V-bounded solution.

\begin{theorem}\label{th:unbs} Let  $\tilde{\Omega}$ be a subset of the domain $\Omega $
and let $$\tilde \Omega^*:=\{(t,x,y)\in \mathbb R\times \mathbb R^{2n}:(t,x)\in \tilde{\Omega},\;(t,y)\in
\tilde{\Omega}\}.$$ Suppose that there exist  functions $V(\cdot):\mathrm{C}^{1}(\mathbb R^{1+n}\!\mapsto
\!\mathbb R)$, $U(\cdot)\!\in\! \mathrm{C}^1(\tilde \Omega^*\!\! \mapsto\! \mathbb R)$, $\eta (\cdot)\!\in\! \mathrm{C}(\mathbb R_+
\!\!\mapsto\! \mathbb R_+)$, and $\beta(\cdot)\!\in\! \mathrm{C}(\mathbb R\times \mathbb R_+ \!\!\mapsto\! \mathbb R_+) $ such
that:

1) the function $V(\cdot)$ is spatially coercive and the function $\eta (\cdot)$ is
positive-definite;

2) the function $\dot U_{(f,f)}(t,x,y):= \frac{\partial U(t,x,y)}{\partial t}+\frac{\partial
U(t,x,y)}{\partial x}\cdot f(t,x)+\frac{\partial U(t,x,y)}{\partial y}\cdot f(t,y)$ satisfies the
inequality
\begin{gather*}
\dot U_{(f,f)}(t,x,y)\ge \beta(t,r)\eta(|U(t,x,y)|)\quad \forall (x,y)\in \tilde
V_t^{-1}((-\infty,r])\cap \tilde \Omega^*_t,
\end{gather*}
with $\tilde  V(t,x,y):=\max\{V(t,x),V(t,y)\}$, and takes positive value at any
point $(t,x,y)\in \tilde \Omega^*$ such that $x\ne y$ and $U(t,x,y)=0$ (if the set of
such points is nonempty);

3) for any sufficiently large $r\ge0$, the functions $\beta(\cdot)$, $h(u):=
\int_{1}^{u}\frac{\mathrm{d}s}{\eta(s)}\;(\;u>0)$, and
  \begin{gather*}
b(t,r):=\max \left\{|U(t,x,y)| :
  \quad  (x,y)\in \tilde  V_t^{-1}((-\infty,r])\cap \tilde \Omega^*\right\}
\end{gather*}
satisfy the conditions
 \begin{gather*}
\int_{0}^{\pm \infty}\beta(s,r)\,\mathrm{d}s=\pm \infty,\quad    \liminf_{t\to \pm
\infty}\frac{h(b(t,r))}{\left|\int_{0}^{t}\beta (s,r)\,\mathrm{d}s\right|}<1.
\end{gather*}

Then the system (\ref{eq:nlsys}) cannot have two different V-bounded solutions
$x(t)$, $y(t)$, $t\in \mathbb R$, whose graphs lie in $\tilde{\Omega}$.
\end{theorem}
\begin{proof} Suppose that the system (\ref{eq:nlsys}) has a pair
of solutions $x(t)$, $y(t),\;t\in \mathbb R$ such that $(t,x(t))$, $(t,y(t))\in
\tilde{\Omega}$ and $x(t)\ne y(t)$ for all $t\in \mathbb R$. Let us show that at least one of
these solutions is not V-bounded.

Using reductio ad absurdum reasoning we suppose that there exists sufficiently
large $r> 0$ such that $|\tilde  V(t,x(t),y(t))|\le r$ for all $t\in \mathbb R$. Consider
the function $u(t):=U(t,x(t),y(t))$.  By condition, the function $u(\cdot)$ is
nondecreasing. Hence, there exist  limits $u_*=\lim_{t\to -\infty}u(t)$,
$u^*=\lim_{t\to \infty}u(t)$ (either finite or infinite).

Firstly, suppose  that  $u_*\ge 0$. If $u(0)=0$, then by condition~2) $\dot
u(0)>0$. Hence, in this case, as well as in the case where $u(0)>0$, we have
the inequality $u(t)>0$ for all $t>0$. Now the condition~2) yields
\begin{gather*}
  h(u(t))-h(u(t_0))\ge \int_{t_0}^{0}\beta(s,r)\,\mathrm{d}s+\int_{0}^{t}\beta(s,r)\,\mathrm{d}s\quad
  \forall t_0>0,\;\forall t\ge t_0.
\end{gather*}
This implies that
\begin{gather*}
   h(b(t,r))-h(u(0))+\int_{0}^{t_0}\beta(s,r)\,\mathrm{d}s\ge \int_{0}^{t}\beta(s,r)\,\mathrm{d}s
   \quad \forall t\ge t_0,
\end{gather*}
and we arrive at contradiction with assumption 3).

Now suppose  that  $u_*<0$. Then there exists  $t'$ such that $u(t')<0$. Thus,
$u(t)\le u(t')$ for all $t<t'$. Then
\begin{gather*}
  \int_{u(t)}^{u(t')}\frac{\mathrm{d}s}{\eta(-s)}\ge \int_{t}^{t'}\beta(s,r)\,\mathrm{d}s\quad \Rightarrow
  \quad h(|u(t)|)-h(|u(t')|)\ge \\\left|\int_{0}^{t}\beta(s,r)\,\mathrm{d}s\right|+\int_{0}^{t'}\beta(s,r)\,\mathrm{d}s
\end{gather*}
from whence, as above, we again arrive at contradiction.
\end{proof}

\begin{remark}\label{rem:uniq2}If $\int_{0}^{1}\frac{1}{\eta(u)}\,\mathrm{d}u<\infty $, then the condition~3)
can be replaced by the following one:
\begin{gather*}
  \liminf_{t\to \infty}\frac{h(b(t,r))+h(b(-t,r))-2h(0)}{\int\limits_{-t}^{t}\beta(s,r) \,\mathrm{d}s}<1.
\end{gather*}
\end{remark}

\section{Studying V-bounded solutions by means of quadratic forms}
\label{QuadForms}

Denote by $\langle \cdot,\cdot \rangle $ a scalar product in $\mathbb R^n$, and let $\|\cdot \|:=\sqrt{\langle
\cdot,\cdot \rangle}$. In this section, the case will be considered where the guiding
function is a time dependent nondegenerate indefinite quadratic form $\langle
S(t)x,x\rangle $. In more detail, the mapping $S(\cdot)\in \mathrm{C}^{1} \left({\mathbb R} \mapsto
\mathrm{Aut}({\mathbb R}^{n})\right)$ assumed to have the following property:
\begin{itemize}
  \item[(a):] for any $t\in \mathbb R$ the operator $S(t)$ is symmetric and there exists a
  decomposition of $\mathbb R^{n}$ into  direct sum of two $S(t)$-invariant subspaces
  $\mathbb L_+(t)$, $\mathbb L_-(t)$ such that the restriction of $S(t)$ on $\mathbb L_+(t)$ (on
  $\mathbb L_-(t)$) is a positive-definite (negative-definite) operator.
  \end{itemize}

Observe that since the subspaces $\mathbb L_+(t),\;\mathbb L_-(t)$ are mutually orthogonal,
the corresponding projectors $P_{\pm}(t):\mathbb R^{n} \mapsto \mathbb L_{\pm}(t)$ are symmetric.

It appears that the function $W(t,x)=\langle S(t)x,x\rangle $ generates a set $\mathcal{W}$
possessing the Wa\.zewski property (C). For the sake of completeness we give
here a proof of the corresponding statement.
\begin{lemma}\label{lem:retrquadform3} Let $W(t,x):=\langle S(t)x,x\rangle $ and let $S(\cdot)$ has the property
(a). Then for any $w>0$, $t_0\in \mathbb R$ there exists a retraction of the set
$W^{-1}(w)$ to the ellipsoid $\{t_0\}\times \left(W_{t_0}^{-1}(w)\cap \mathbb L_+(t_0)\right)$.
 \end{lemma}

\begin{proof}
From $S(t)$-invariance of subspaces $\mathbb L_+(t),\;\mathbb L_-(t)$ it follows that
$P_{\pm}(t)S(t)=S(t)P_{\pm}(t)$ and, as a consequence, we have the representation
\begin{gather*}
  S(t)=(P_+(t)+P_-(t))S(t)(P_+(t) +P_-(t))=\\P_+(t)S(t)P_+(t) +P_-(t)S(t)P_-(t).
\end{gather*}
Put
\[
  S_+(t):=P_+(t)S(t)P_+(t),\;S_-(t)=:P_-(t)S(t)P_-(t)
\]
Obviously, the kernel of the operator $S_+(t)$ (operator $S_-(t)$) is the
subspace  $\mathbb L_-(t)$ (subspace $\mathbb L_+(t)$), and the restriction of this operator
on $\mathbb L_+(t)$ (on $\mathbb L_-(t)$) is a positively definite (negatively definite)
operator.

Now observe that for arbitrary  $t\in \mathbb R$ and $w>0$ there exists a retraction
of $W_{t}^{-1}(w)=\{x\in \mathbb R^n:\langle S(t)x,x\rangle =w\}$ to the intersection of this
set with the subspace $\mathbb L_+(t)$. In fact, one can define such a retraction
 by a mapping $x \mapsto \varpi(t,x)P_+(t)x$, provided that the scalar function
 $\varpi(t,x)$ is determined from condition $\langle S_+(t)\varpi(t,x)x,\varpi(t,x)x\rangle =w$ for
all $x\in W_{t}^{-1}(w)$. Since $w>0$, then $W_{t}^{-1}(w)\cap
\mathbb L_{-}(t)=\varnothing$, and hence, $\langle S_+(t)x,x\rangle >0$ for all $x\in
W_{t}^{-1}(w)$. Therefore
\begin{gather*}
  \varpi(t,x)=\sqrt{\frac{w}{\langle S_+(t)x,x\rangle}}.
\end{gather*}

Now it remains only to show that the set  $\{t_0\}\times
W_{t}^{-1}(w)=W^{-1}(w)\cap\Pi_{t_0}$ is a retract of $W^{-1}(w)$. Introduce the
operator $R(t):=\sqrt{S^2(t)}=S_+(t) -S_-(t)$. Then we get
\begin{gather*}
  S(t)=R(t)(P_+(t) -P_-(t))=(P_+(t) -P_-(t))R(t).
\end{gather*}
The quadratic form  $\langle S(t)x,x\rangle $ by means of the substitution
$x=\left[\sqrt{R(t)}\right]^{-1}\!\!\!y$ is reduced to  $\langle (P_+(t)-P_-(t))y,y\rangle$.
Obviously, $P_+(t)-P_-(t)$ is  a symmetric orthogonal inversion operator:
$$(P_+(t)-P_-(t))^*=P_+(t)-P_-(t),\quad  (P_+(t)-P_-(t))^2=E.$$

 From the representation  via the Riesz formula (see, e.g.,
\cite[c.~34]{DalKre70}) it follows that the projectors $P_{\pm}(t)$  smoothly
depend on parameter. Therefore the mutually orthogonal subspaces $\mathbb L_+(t)$ and
$\mathbb L_-(t)$ have constant dimensions $n_+$,  $n_{-}$ and define smooth curves
$\gamma_+$, $\gamma_-$ in Grassmannian  manifolds $G(n,n_{+})$ and $G(n,n_{-})$
respectively. Since $G(n,n_{+})$ is a base space of a principal fiber bundle,
namely, $G(n,n_{+})=O(n)/O(n_+)\times O(n_-)$, then there exists a smooth curve
$Q(t)$ in $O(n)$, which is projected onto  $\gamma_+(t)$, the operator $Q(t_0)$
being the identity element $E$  of the group $O(n)$. Obviously,
$\mathbb L_+(t)=Q(t)\mathbb L_+(t_0)$ and, as a consequence,
\begin{gather*}
  P_{\pm}(t)=Q(t)P_{\pm}(t_0)Q^{-1}(t).
\end{gather*}

From the above reasoning it follows that the change of variables
\begin{gather*}
  x=\left[\sqrt{R(t)}\right]^{-1}Q(t)\sqrt{R(t_0)}y
\end{gather*}
reduces the quadratic form  $W(t,x):=\langle S(t)x,x\rangle $ to $W(t_0,y)=\langle S(t_0)y,y\rangle
$, and then the mapping
\begin{gather*}
  \mathbb R\times \mathbb R^{n} \mapsto\{t_0\}\times \mathbb R^n:\quad  (t,x)  \mapsto
  \left(t_0,\sqrt{R(t)}Q^{-1}(t)\left[\sqrt{R(t_0)}\right]^{-1}x\right)
\end{gather*}
determines a retraction of the set $W^{-1}(w)$ to the set $W^{-1}(w)\cap
\Pi_{t_0}$.
\end{proof}

Now consider the quasilinear system
\begin{equation}\label{eq:quaslin}
  \dot{x}=f(t,x):=A(t)x+g(t,x)
\end{equation}
and assume that the following conditions hold:
\begin{itemize}
  \item[(b):] the mapping $A(\cdot)\in \mathrm{C}\left({\mathbb R} \mapsto {\mathrm{Hom}}\left({\mathbb R}^{n} \right)\right)$ is
   such that $\sup_{t\in \mathbb R}\|A(t)\|=:a <\infty$  and the linear system $\dot{x}=A(t)x$ is
   exponentially dichotomic on  $\mathbb R$;  i.e. there exists a  mapping $C(\cdot)\in \mathrm{C}^{1} \left({\mathbb R} \mapsto \mathrm{Aut}({\mathbb R}^{n})\right)$
  possessing the property (a) with $S(t)=C(t)$, and, in addition,
\begin{gather*}
  \sup_{t\in \mathbb R}\|C(t)\|=:c <\infty,\quad \inf_{t\in \mathbb R}|\det C(t)|=:\sigma >0, \\
  \left\langle(2C(t)A(t)+\dot C(t))x,x\right\rangle
\ge \|x\|^2\quad \forall t\in \mathbb R,\;x\in \mathbb R^{n}
\end{gather*}
(see, e.g. \cite{MSK03,Sam02});
  \item[(c):] there exist $k >0$  and  $\varphi (\cdot)\in
\mathrm{C}^1({\mathbb R} \mapsto (0,\infty))$ such that $\sup_{t\in \mathbb R}\frac{|\dot \varphi(t)|}{\varphi(t)}=:l
<\infty $ and the mapping $ g(\cdot)\in \mathrm{C}\left({\mathbb R}^{1+n} \mapsto {\mathbb R}^{n} \right)$ satisfies
the inequality  $$ \left\| g(t,x)\right\| \le k\left\| x\right\| +\varphi (t)\quad \forall
(t,x)\in {\mathbb R}^{1+n}. $$
\end{itemize}

The well known approach to establish sufficient conditions for the existence of
bounded solutions to \eqref{eq:quaslin} is based on the method  of integral
equations which allows to apply different versions of fixed point theorems
(see, e.g. \cite{MNS72,Per03}). Our goal is to show that by means of V-W-pair
one can not only establish the existence of bounded solutions (in the case
where $\varphi(t)$ is bounded), but also show how their asymptotic behavior depends
on $\varphi(t)$ as $t\to \pm \infty $.

For any $t\in \mathbb R$, put
\begin{gather*}
  \lambda^+_C(t):=\max_{\|x\|=1} \langle C(t)x,x\rangle,\quad
  \lambda^-_C(t):=\min_{\|x\|=1} \langle C(t)x,x\rangle,\\
    \lambda^+_{C,\mathrm{min}}(t):=\min \left\{\langle C(t)x,x\rangle:\|x\|=1,x\in \mathbb L_+(t)\right\}
  \end{gather*}
and
\begin{gather*}
  F(r):=\begin{cases}
    \frac{d}{m}r^2+2\left( \frac{c}{m}+\frac{d}{m^2}\right)
    \left(\frac{1}{m}\ln(mr+1)-r\right) &\text{if}\quad
    r\ge \frac{c}{d}, \\
    F \left(\frac{c}{d}\right) & \text{if}\quad 0\le r <\frac{c}{d},
  \end{cases}
\end{gather*}
where $d:=\frac{1}{2}-c(k+l)$, $m:=a+k+l$.

\begin{theorem}\label{th:bsqlsys}
Let  the conditions~(b),(c) hold true and let the numbers $c,k,l$ satisfy the
inequality
\begin{gather*}
  c(k+l)<\frac{1}{2}.
\end{gather*}
Then the system \eqref{eq:quaslin} has a solution $x_*(t),\;t\in \mathbb R$, such that
\begin{gather}\label{eq:estbsqls}
  \|x_*(t)\|\le r_*\varphi(t)\quad \forall t\in \mathbb R
\end{gather}
where $r_*$ is the root of equation
\begin{gather*}
  F(r)= F\left(\frac{c}{d}\right)+
  \frac{c^2}{2d^2}\sup_{t\in \mathbb R}\left[\sup_{s\ge t}\lambda^+_C(s)-\inf_{s\le t}\lambda^-_C(s)\right].
\end{gather*}
If, in addition,
\begin{gather*}
  \|g(t,x)-g(t,y)\|\le k\|x-y\|\quad \forall (t,x,y)\in \mathbb R^{1+2n},
\end{gather*}
then $x_*(t)$ is a unique solution of the system \eqref{eq:quaslin} for which
the ratio $\frac{\|x\|}{\varphi(t)}$ is bounded on $\mathbb R$.
\end{theorem}
\begin{proof}
First, we show that the system \eqref{eq:quaslin} has the following V-W-pair
\begin{gather}\label{eq:VWfi}
  V(t,x):= F\left(\frac{\|x\|}{\varphi(t)}\right)-F\left(r_0\right),\quad
  \quad W(t,x)=\frac{\langle C(t)x,x\rangle}{\varphi^2(t)}
\end{gather}
where $r_0$ is an arbitrary number greater than ${c}/{d}$. In fact, from the
inequalities \begin{gather*}\left|\frac{\mathrm{d}}{\mathrm{d}t} \left[\frac{\left\|
x\right\|^{2}}{\varphi^{2} (t)} \right]_{\!\!f} \right|\le \frac{2}{\varphi^{2} (t)} \left[(a+k)\left\|
x\right\|^{2}+\varphi(t)\|x\| \right]+2\frac{|\dot{\varphi}(t)|}{\varphi^3(t)}\|x\|^2 \le \\
2m\frac{\left\| x\right\|^{2}}{\varphi^{2} (t)}+ 2\frac{\left\| x\right\|}{\varphi (t)},\end{gather*}
$$\begin{array}{c} \dfrac{\mathrm{d}}{\mathrm{d}t}\left[\dfrac{\langle C(t)x,x\rangle}{\varphi^{2} (t)}
\right]_{\!\!f}   \ge \dfrac{1}{\varphi^{2} (t)} \left[(1-2ck)\| x \|^{2} -2c\varphi (t)\|
x\|\right]-2c \dfrac{|\dot{\varphi}(t)|}{\varphi^{3} (t)} \| x\|^{2} \ge\\
2d\dfrac{\|x\|^2}{\varphi^2(t)}-2c\dfrac{\|x\|}{\varphi(t)}
\end{array}$$
and equality
\begin{gather*}
  F(r)-F(r_0)=2\int_{r_0}^{r}\frac{ds^2-cs}{ms+1}\,\mathrm{d}s,
\end{gather*}
it follows  that $\dot W_f(t,x)> 2(dr^2_0-cr_0)>0$ and $|\dot V_f(t,x)|\le \dot
W_f(t,x)$ once $\|x/\varphi(t)\|> r_0>c/d$, or, equivalently, $V(t,x)>0$.

Next, it is easily seen that in our case
\begin{gather*}
  w^0(t)=r^2_0\lambda^+_C(t),\quad w_0(t)=r^2_0\lambda^-_C(t).
\end{gather*}
If we pick $w_*,w^*$ in such a way that
\begin{gather*}
  w_*<r^2_0\inf_{t\in \mathbb R} \lambda^-_C(t),\quad w^*>r^2_0\sup_{t\in \mathbb R} \lambda^+_C(t),
\end{gather*}
then, in view of Lemma~\ref{lem:retrquadform3}, to satisfy the
conditions~(A),(B),(C) it is sufficient to define
\begin{gather*}
  \mathcal{W}:=W^{-1}\left((w_*,w^*)\right),\quad \mathcal{M}_t:=W_t^{-1}\left([0,w^*]\right)\cap \mathbb L_+(t).
\end{gather*}
Note, that in our case $c_*=c^*=1$ and $\alpha(t)\ge 2(dr^2_0-cr_0)>0$.

Lastly, from (b) it follows that $\inf_{t\in \mathbb R} \lambda^+_{C,\mathrm{min}}(t):=\sigma_1>0$.
Hence, $\frac{\|x\|^2}{\varphi^2(t)}\le \frac{w^*}{\sigma_1}$ for all $t\in \mathbb R$, all
$x\in \mathcal{M}_t $, and this yields $\nu <\infty $. Now, by the Theorem~\ref{th:exbs},
there exists a solution $x_*(t),\;t\in \mathbb R$, of the system~\eqref{eq:quaslin}
such that
\begin{gather*}
  V(t,x_*(t))\le \frac{r^2_0}{2} \left[\sup_{s\ge t}\lambda^+_C(s)-\inf_{s\le
  t}\lambda^-_C(s)\right].
\end{gather*}
The estimate \eqref{eq:estbsqls} is easily obtained by letting $r_0$ tend to
$c/d$.

In order to prove the uniqueness of $x_*(t)$, it remains only to apply the
Theorem~\ref{th:unbs} in the case where $U(t,x,y):=W(t,x-y)$,
$V=V_+(t,x):=\|x\|^2/\varphi^2(t)$, $\eta(u):= u$, $\beta(t,r)=:\left(1-2(k+l)\right)/c$,
$b(t,r):=4cr$.
\end{proof}

\begin{remark} The number $r_*$ does not exceed the largest root of the quadratic equation
\begin{gather*}
  \frac{d}{m}r^2-2\left( \frac{c}{m}+\frac{d}{m^2}\right)r-F\left(\frac{c}{d}\right)
  -\frac{c^2}{2d^2}\sup_{t\in \mathbb R}\left[\sup_{s\ge t}\lambda^+_C(s)-\inf_{s\le
  t}\lambda^-_C(s)\right]=0.
\end{gather*}
\end{remark}
\begin{remark}\label{rem:V*}The assertion of the Theorem~\ref{th:bsqlsys} remains
true if we require that the function $g(\cdot)$ is defined and satisfies the
condition~(c) not on the whole $\mathbb R^{1+n}$ but only on a domain $\Omega $ which
contains the set $W^{-1}\left([w_*,w^*]\right)\cap V^{-1}([0,V^*])$ where V--W-pair is
defined by \eqref{eq:VWfi} and $V^*=w^*+\max\{\nu,-w_*\}$.
\end{remark}

\begin{example}\label{ex:vbs1}
Consider the following singular boundary value problem for scalar second order
differen\-tial equation
\begin{gather}\label{eq:sbvp1}
  \frac{\mathrm{d}}{\mathrm{d}t}\left(\frac{\dot z}{\rho(t)}\right)-\omega (t)z=Z (t,z,\dot z),\\ \label{eq:sbvp2}
  z(-\infty)=z(+\infty)=0,
\end{gather}
where $\rho(\cdot)\in \mathrm{C}^{1}(\mathbb R\!\mapsto \!(0,\infty))$, $\omega (\cdot)\in \mathrm{C}(\mathbb R\!\mapsto
\!(0,\infty))$ are bounded functions and the function $Z (\cdot)\in \mathrm{C}(\mathbb R^{3}\!\mapsto
\!\mathbb R)$ satisfies a global Lipschitz condition: there exists a constant $\ell$
such that
\begin{gather*}
|Z(t,x_1,y_1)-Z(t,x_2,y_2)|\le
  \ell\sqrt{(x_1-x_2)^2+(y_1-y_2)^2}\\ \forall \{t,x_1,y_1,x_2,y_2\}\subset \mathbb R.
\end{gather*}

 Let us show that if there exists a function $\varphi(\cdot)\in \mathrm{C}^{1}(\mathbb R\!\mapsto
\!(0,\infty))$ such that
\begin{gather*}
\quad|Z(t,0,0)|\le \varphi(t),\quad\lim_{|t|\to \infty}\varphi(t)=0,\quad \sup_{t\in
\mathbb R}\frac{|\dot \varphi(t)|}{\varphi(t)}:= l<\infty
  \end{gather*}
  and
\begin{gather*}
  k+l<\delta
\end{gather*}
where
\begin{gather*}
\delta := \min\left\{\inf_{t\in \mathbb R}\rho(t),\inf_{t\in \mathbb R}\omega (t)\right\},\quad k:=\ell\max
\left\{1,\sup_{t\in \mathbb R}\rho(t)\right\},
\end{gather*}
then the problem~\eqref{eq:sbvp1}--\eqref{eq:sbvp2} has a unique solution
$z_*(t)=O(\varphi(t))$.

By letting  $x_1=z,\;x_2=\dot z/\rho(t)$, the equation~\eqref{eq:sbvp1} becomes
equivalent to 2-D system of the form \eqref{eq:quaslin}  with
\begin{gather*}
  A(t)\!=\!
  \begin{pmatrix}
    0 & \rho(t) \\
    \omega(t) & 0
  \end{pmatrix},\; g(t,x)\!=\!
  \begin{pmatrix}
    0 \\
         Z (t,x_1,\rho(t)x_2)
  \end{pmatrix}
\end{gather*}

Set $\langle C(t)x,x\rangle =\frac{x_1x_2}{\delta}$. Obviously this is a nondegenerate
indefinite quadratic form of Morse index 1. One can easily show that
$\|C(t)\|=\frac{1}{2\delta}=:c$, $\frac{\mathrm{d}}{\mathrm{d}t}\langle C(t)x,x\rangle_{A(t)x}\ge \|x\|^2
$, $c(k+l)<1/2$, $\|g(t,0)\|\le \varphi(t)$, and $\|g(t,x)-g(t,y)\|\le k\|x-y\|$.
Now the unique solvability of the problem~\eqref{eq:sbvp1}--\eqref{eq:sbvp2} in
the class of functions $z(t)=O(\varphi(t))$ follows from the
Theorem~\ref{th:bsqlsys}.

Note that if we slightly simplify our task by replacing the
condition~\eqref{eq:sbvp2} with $\sup_{t\in \mathbb R}|z(t)|<\infty $, then the
sufficient condition for solvability of the corresponding problem takes the
form
\begin{gather*}
  \sup_{t\in \mathbb R}|Z(t,0,0)|<\infty,\quad k<\delta
\end{gather*}
(obviously, in this case $\varphi(t)\equiv \mathrm{const}$, and $l=0$). At the same time,
by applying results of \cite{Per03} combined with estimates for Green function
derived in \cite{MSK03,Sam02}, we can only obtain a rougher condition
\begin{gather*}
  2k\delta^{-3/2}\sqrt{\max \left\{\sup_{t\in \mathbb R}\rho(t),\sup_{t\in \mathbb R}\omega(t)\right\}}<1
\end{gather*}
(note that the expression under the square root is not less than $\delta $).
\end{example}

Now let us lay down sufficient conditions for the existence of $\mathrm{V}_+$-bounded
solutions in the case where $f(\cdot)\in \mathrm{C}(\mathbb R^{1+n} \mapsto \mathbb R^{n})$ is essentially
nonlinear, e.g., $\|f(t,x)\|/\|x\|\to \infty,\;x\to \infty$. We are going to construct
a V-W-pair in the form $V(t,x)=F(V_+(t,x))$, $W(t,x)=\langle S(t)x,x\rangle $,
$V_+(t,x)=\langle B(t)x,x\rangle $ under the following conditions:
\medskip
\begin{itemize}
\item[(d):] for any $t\in \mathbb R$, the operator $B(t)$ is positively
definite and  there exist projectors $P_+(t), P_-(t)$ on invariant subspaces
$\mathbb L_{+}(t),\mathbb L_{-}(t)$ of operator $S(t)$ such that the restriction of $S(t)$ on
$\mathbb L_+(t)$ (on $\mathbb L_-(t)$) is a positively definite (negatively definite)
operator.
\item[(e):] there exist functions $\gamma(\cdot)\in \mathrm{C}(\mathbb R\!\mapsto \!(0,\infty))$, $\Gamma (\cdot)\in
  \mathrm{C}((0,\infty)\!
  \mapsto\! \mathbb R)$, $\Delta (\cdot)\in \mathrm{C}((0,\infty)\! \mapsto\!(0,\infty))$ such that
  \begin{gather*}
  \min_{\{ x\in \mathbb R^{n} :\langle B(t)x,x\rangle =v\}}\langle S(t)f(t,x),x\rangle\ge \gamma(t)\Gamma(v)\quad \forall v>0,\\
  \max_{\{x\in \mathbb R^{n} :\langle B(t)x,x\rangle =v\}}\left|\langle B(t)f(t,x),x\rangle\right|\le \gamma(t)\Delta(v)\quad \forall v>0,
\end{gather*}
and
\begin{gather*}
  \int_{-\infty}^{0}\!\!\gamma
(t)\,\mathrm{d}t=\int_{0}^{\infty}\!\!\gamma (t)\,\mathrm{d}t=\infty;
\end{gather*}
\item[(f):] the following inequalities hold true
\begin{gather*}
  \sup_{t\in \mathbb R}\lambda^+(t)<\infty,\quad \inf_{t\in \mathbb R}\lambda_{-}(t)>-\infty,
  \quad \limsup_{t\to -\infty}\lambda_-^+(t)>0,\\
  \inf_{t\in \mathbb R}\frac{\mu_-(t)}{\gamma (t)}=:\xi >-\infty,\quad \sup_{t\in
  \mathbb R}\frac{M(t)}{\gamma (t)}=:\varsigma <\infty
\end{gather*}
where $\lambda^+(t)$ and $\lambda_-(t)$ are, respectively, the maximal and the minimal
characteristic values of the pencil $S(t)-\lambda B(t)$, $\lambda_-^+(t)$ is the minimal
characteristic value of the pencil $ P_+(t)\left[S(t)-\lambda
B(t)\right]\left|_{\mathbb L_+(t)}\right.$, $M(t)$ is the maximum of moduli of maximal and
minimal characteristic values of the pencil $\dot{B}(t)-\mu B(t)$, and  $\mu_-(t)$
is the minimal characteristic value of the pencil $\dot{S}(t)-\mu B(t)$.
\item[(g):] there exists a number $v_0>0$ such that
\begin{gather*}
  2\Gamma(v_0)+\xi v_0>0,\quad\frac{\Gamma(v)-\Gamma(v_0)}{v-v_0}\ge -\frac{\xi}{2}\quad \forall
  v>v_0,\\
   \int_{v_0}^{\infty}\frac{2\Gamma(v)+\xi v}{2\Delta(v)+\varsigma v}\,\mathrm{d}v=\infty.
\end{gather*}
 \end{itemize}

\medskip

We arrive at the following result.

\begin{theorem}\label{th:thbs1}
Let  the system   (\ref{eq:nlsys}) satisfies   in $\Omega :=\mathbb R^{1+n}$ the
conditions~(d)--(g).  Then there exists a solution $x_*(t)$ of this system such
that
\begin{gather*}
\langle B(t)x_*(t),x_*(t)\rangle  \le
  F^{-1}\left(\frac{v_0}{2}\left[\sup_{s\ge t}\lambda^+(s)-\inf_{s\le t}\lambda_-(s)\right]\right)\le v_*\quad
 \forall   t\in \mathbb R
\end{gather*}
where
\begin{gather*}
  F(v):=\int_{v_0}^{v}\frac{2\Gamma(u)+\xi
u}{2\Delta(u)+\varsigma u}\,\mathrm{d}u
\end{gather*}
and $v_*$ is the root of the equation $$F(v)=\frac{v_0}{2}\left[\sup_{t\in
\mathbb R}\lambda^+(t)-\inf_{t\in \mathbb R}\lambda_-(t)\right].$$

If in addition $ 2\gamma(t)+\mu_-(t)>0$ for all $t\in \mathbb R$,
\begin{gather*}
  \left\langle S(t)\left(f(t,x+y)-f(t,x)\right),y\right\rangle \ge
  \gamma(t)\langle B(t)y,y\rangle\quad \forall (t, x,y)\in \mathbb R^{1+2n},
\end{gather*}
and
\begin{gather*}
\int_{0}^{\pm \infty}\frac{2\gamma(s)+\mu_-(s)}{\max\{\lambda^+(s),|\lambda_-(t)|\}}\,\mathrm{d}s=\pm \infty,\quad
\liminf_{t\to \pm \infty}\frac{\ln\max\{\lambda^+(t),|\lambda_-(t)|\}}{
\left|\int_{0}^{t}\frac{2\gamma(s)+\mu_-(s)}{\max\{\lambda^+(s),|\lambda_-(t)|\}}\,\mathrm{d}s\right|}<1,
\end{gather*}
then $x_*(t)$ is a unique solution of the system (\ref{eq:nlsys}) satisfying
the condition $$\sup_{t\in \mathbb R}\langle B(t)x_*(t),x_*(t)\rangle <\infty.$$
\end{theorem}

\begin{proof} Put $W(t,x):=\langle S(t)x,x\rangle $, $V_+(t,x):=\langle B(t)x,x\rangle$. Since
\begin{gather*}
  M(t)=\max_{\{x\in \mathbb R^{n}:V_+(t,x) =1\}}\left|\langle \dot{B}(t)x,x\rangle\right|,\quad
  \mu_-(t)=\min_{\{x\in \mathbb R^{n}:V_+(t,x) =1\}}\langle \dot{S}(t)x,x\rangle
\end{gather*}
(see, e.g., \cite{Gan67}), then
\begin{gather*}
 \left|\left[\dot{V}_+(t,x)\right]_{\!f}\right|\le 2\gamma(t)\Delta(V_+(t,x))+M(t)V_+(t,x)\le \\
 \gamma(t)\left(2\Delta(V_+(t,x))+\varsigma V_+(t,x)\right),\\
   \dot{W}_f(t,x)\ge 2\gamma(t)\Gamma(V_+(t,x))+\mu_-(t)V_+(t,x)\ge \\
   \gamma(t)\left(2\Gamma(V_+(t,x))+\xi V_+(t,x)\right),
\end{gather*}
once $V_+(t,x) >v_0$, and it is naturally to define in this case
\begin{gather}\label{eq:gGa}
  V(t,x)=F(V_+(t,x)).
\end{gather}
Obviously that the inequality~\eqref{eq:mainineq}  and condition~(B) are
satisfied with $c_*=c^*=1$, $\alpha(t)\ge \gamma(t)(2\Gamma(v_0)+\xi v_0)$.

 Taking into account that the function $W_t(x)$
has the unique critical point $x=0$, we have
\begin{gather*}
  w^0(t):=\max_{\{x\in \mathbb R^{n} :V_+(t,x) \le
  v_0\}}W(t,x)=\lambda^+(t)v_0,\\
w_0(t):=\min_{\{x\in \mathbb R^{n} :V_+(t,x) \le v_0\}}W(t,x)=\lambda_-(t)v_0.
\end{gather*}
If we choose numbers $w_*,w^*$ in such a way that
\begin{equation}\label{eq:w_+}
  w_*<\omega_0=\inf_{t\in \mathbb R}\lambda_-(t)v_0,\quad w^*>\omega^0=\sup_{t\in \mathbb R}\lambda^+(t)v_0,
\end{equation}
and define $ \mathcal{W}:=W^{-1}\left((w_*,w^*)\right)$, then the condition~(A) will be
satisfied.

As has been already shown in proof of Theorem\ref{th:bsqlsys} the family of
sets $\mathcal{M}_t:=W_t^{-1}\left([0,w^*]\right)\cap \mathbb L_+(t)$ satisfy the condition~(C). Now
to prove the existence of V-bounded solution it remains only to show that $\nu
<\infty $. It is easily seen that
\begin{gather*}
  \min\{W_{t}(x):x\in \mathcal{M}_t,\; V_+(t,x) >v_0\}=\lambda_-^+(t)v_0>0,\\
  \max_{\{x\in \mathcal{M}_t\}}V_+(t,x)
  =\frac{w^*}{\lambda_-^+(t)},
\end{gather*}
and in view of condition~(f) we have $\liminf_{t\to -\infty}(w^*/\lambda_-^+(t))<\infty$.
Hence, $\nu <\infty $.

I order to prove the uniqueness of V$_+$-bounded solution of the system
\eqref{eq:nlsys}, introduce the function $U(t,x,y):=\langle S(t)(x-y),(x-y)\rangle $. It
is easily seen that
\begin{gather*}
  \dot U_{(f,f)}(t,x,y)\ge (2\gamma(t)+\mu_-(t))\langle B(t)(x-y),x-y\rangle\ge\\
  \frac{2\gamma(t)+\mu_-(t)}{\max\{\lambda^+(t),|\lambda_-(t)|\}}|U(t,x,y)|,\\
|U(t,x,y)|\le \max\{\lambda^+(t),|\lambda_-(t)|\}\langle B(t)(x-y),x-y\rangle.
\end{gather*}
Now the uniqueness result follows from Theorem~\ref{th:unbs} if we define
\begin{gather*}
  \beta(t,r):=\frac{2\gamma(t)+\mu_-(t)}{\max\{\lambda^+(t),|\lambda_-(t)|\}},\\ b(t,r):=
  4\max\{\lambda^+(t),|\lambda_-(t)|\}r, \quad \eta(u):=u.
\end{gather*}
\end{proof}

\section{V-bounded solutions of Lagrangian systems}\label{Lagrsys}

Consider a natural Lagrangian system  subjected to smooth time-varying
holonomic constraint. The Lagrangian of such a system can be represented in the
form
\begin{equation}\label{eq:Ltqq}
  L(t,q,\dot q):=\tfrac{1}{2}\langle \mathcal{A}(t,q)\dot q,\dot q\rangle+ \langle a(t,q),\dot q\rangle +\Phi(t,q)
\end{equation}
where $q=(q_1,\ldots,q_m)\in \mathbb R^{m}$ are generalized coordinates, $\mathcal{A}(\cdot):\mathbb R^{1+m}
\mapsto \mathrm{Aut}(\mathbb R^{m})$, $a(\cdot):\mathbb R^{1+m} \mapsto \mathbb R^{m}$, $\Phi(\cdot):\mathbb R^{1+m} \mapsto \mathbb R$ are
$\mathrm{C}^{2}$-mappings, and besides, $\mathcal{A}(\cdot)$ takes values in the space of
positive-definite operators. Our goal is to show that if the Lagrangian has
certain directional quasiconvexity property, namely

\medskip

\begin{itemize}
  \item[($\alpha $):] there exist positive numbers $\kappa $, $R$ and a  function
$\Psi(\cdot):\mathbb R^{1+m} \mapsto \mathbb R_+$ such that from $$\tfrac{1}{2}\langle \mathcal{A}(t,q)\dot q,\dot q\rangle+
\Psi (t,q)\ge R$$ it follows that
\begin{gather}\label{eq:quasiconv}
  \hspace*{4mm} \dfrac{\partial L}{\partial q_i}q_i+\dfrac{\partial L}{\partial \dot q_i}\dot q_i \ge \kappa
\left(\tfrac{1}{2}\langle \mathcal{A}(t,q)\dot q,\dot q \rangle +\Psi (t,q)\right)\\(\text{summation over
repeating indices}),\nonumber
\end{gather}
\end{itemize}
\medskip

\noindent then under some additional technical growth conditions imposed on
$\mathcal{A}(\cdot),a(\cdot),\Psi(\cdot)$ the Lagrangian system possesses a global solution along
which the function $\frac{1}{2}\langle \mathcal{A}(t,q)\dot q,\dot q\rangle+\Psi (t,q)$ is bounded.

\begin{remark}\label{rem:quasiconv} It is easilily seen that the
inequality \eqref{eq:quasiconv} yields
\begin{gather} \label{eq:quasiconvA}
  \left\langle \left(\mathcal{A}(t,q)+\frac{1}{2}\frac{\partial \mathcal{A}(t,q)}{\partial q_i}q_i\right)y,y\right\rangle \ge
  \frac{\kappa}{4}\langle \mathcal{A}(t,q)y,y \rangle \quad \forall(t,q,y)\in \mathbb R^{1+2n}.
\end{gather}
\end{remark}

It should be also noted that the Assumptions (H4),(H5) in \cite{Cie03b} implies
that
\begin{gather*}
  \dfrac{\partial L}{\partial q_i}q_i+\dfrac{\partial L}{\partial \dot q_i}\dot q_i \ge \kappa
\left(\|\dot q\|^2 +\|q\|^2\right)
\end{gather*}
once $\|\dot q\|^2+\|q\|^2$ is sufficiently large.
\medskip

In what follows, we shall also assume that:
\medskip

\begin{itemize}

\item[($\beta $):] there exists a nondecreasing coercive functions
$\underline{\Theta} (\cdot):\mathbb R_+ \mapsto \mathbb R_+$, $\overline{\Theta} (\cdot):\mathbb R_+ \mapsto \mathbb R_+$ such
\begin{gather*}
 \underline{\Theta} \left(\Psi (t,q)\right) \le \langle \mathcal{A}(t,q)q,q\rangle \le \overline{\Theta} \left(\Psi (t,q)\right)
  \quad \forall(t,q)\in \mathbb R^{1+m}.
\end{gather*}
\item[($\gamma $):] there exist numbers $\theta \in [0,1]$ and $K>0$
such that from $$\tfrac{1}{2}\langle \mathcal{A}(t,q)\dot q,\dot q\rangle+ \Psi (t,q)\ge R$$ it
follows that
\begin{gather*}
\left|  \frac{1}{2}\left\langle \frac{\partial \mathcal{A}(t,q)}{\partial t}\dot q,\dot q\right\rangle+
  \frac{\partial (\Phi (t,q)+\Psi(t,q))}{\partial q_{i}}\dot q_i+\frac{\partial \Psi (t,q)}{\partial t}\right|\le\\
  K\left(\frac{1}{2}\langle \mathcal{A}(t,q)\dot q,\dot q\rangle+ \Psi (t,q)\right)^{\theta +1},
\end{gather*}
$R$ being the same number as in ($\alpha $).

\item[($\delta $):] there exist a
nondecreasing function  $\Xi (\cdot):\mathbb R \mapsto\mathbb R_+$ such that
 \begin{gather*}
 \max_{\|y\|=1}\frac{\left|\langle a(t,q),y\rangle\right|}{\sqrt{\langle
\mathcal{A}(t,q)y,y\rangle}}\le \Xi (\Psi (t,q))\quad  \forall (t,q)\in \mathbb R^{1+m}.
\end{gather*}
\end{itemize}

In order to apply the results of Section~\ref{ExistUniqVBS}, introduce the
functions
\begin{equation}\label{eq:defofVW}
\begin{split}
  W(t,q,\dot q) & :=\left\langle \mathcal{A}(t,q)\dot q+a(t,q),q\right\rangle, \\
  V(t,q,\dot q) &  :=\bar V \left(\tfrac{1}{2}\langle \mathcal{A}(t,q)\dot q,\dot
  q\rangle +\Psi (t,q)\right)
\end{split}
\end{equation}
where $\bar V(\cdot)\in \mathrm{C}^1(\mathbb R \mapsto(-1,\infty))$ is a strictly increasing function
which for $r\ge R$ is defined as
\begin{gather*}
  \bar V(r)=:
\begin{cases}
    \ln(r/R) &\; \text{if}\quad \theta =1, \\
    \left(r^{1-\theta}-R^{1-\theta}\right)/(1-\theta) & \;\text{if}\quad \theta \in [0,1).
  \end{cases}\end{gather*}

\begin{lemma}\label{lem:VWpair}  From $V(t,q,\dot q)\ge0$ it follows that
\begin{gather*}
    \left|\dot V_f(t,q,\dot q)\right|\le \frac{K}{\kappa}\dot W_f(t,q,\dot q)\quad\text{and}\quad
  \dot W_f(t,q,\dot q)\ge  \kappa R
\end{gather*}
where $f(t,q,\dot q):= \left(\!\dot q,\left(\tfrac{\partial^2L}{\partial \dot
q^2}\right)^{-1}\left(\tfrac{\partial L}{\partial q}-\tfrac{\partial^2L}{\partial t\partial \dot
q}-\tfrac{\partial^2L}{\partial \dot q\partial q_i}\dot q_i\right)\!\!\right)$ is the vector field
generated in the phase space $\mathbb R^{2m}$ by the Lagrangian system.
\end{lemma}
\begin{proof}
Note that $W= \tfrac{\partial L}{\partial \dot q_i}q_i $. The equation of motion
\begin{gather*}
  \frac{\mathrm{d}}{\mathrm{d}t}\dfrac{\partial L}{\partial \dot q}=\dfrac{\partial L}{\partial q}
\end{gather*}
yields
\begin{gather*}
  \frac{\mathrm{d}}{\mathrm{d}t} \dfrac{\partial L}{\partial \dot q_i}q_i\left|_{\!f}\right. =\dfrac{\partial L}{\partial q_i}q_i+
  \dfrac{\partial L}{\partial \dot q_i}\dot q_i.
\end{gather*}
Obviously,
\begin{gather*}
  V(t,q,\dot q)\ge 0\quad \Leftrightarrow\quad
  \tfrac{1}{2}\langle \mathcal{A}(t,q)\dot q,\dot q \rangle +\Psi (t,q)\ge R.
\end{gather*}
Then by assumption~$(\alpha)$ we have
\begin{gather}\label{eq:WgeR} \dot W_f(t,q,\dot q)\ge \kappa
\left(\tfrac{1}{2}\langle \mathcal{A}(t,q)\dot q,\dot q \rangle +\Psi (t,q)\right)\ge \kappa R\quad
\text{once}\quad V(t,q,\dot q)\ge 0.
\end{gather}

In order to estimate the function $\dot V_f(\cdot)$, introduce the Hamiltonian in a
standard way:
\begin{gather*}
  H(t,q,\dot q)=\dfrac{\partial L}{\partial \dot q_i}\dot q_i-L=\frac{1}{2}\langle \mathcal{A}(t,q)\dot q,\dot
  q\rangle -\Phi(t,q).
\end{gather*}
As is well known, $\frac{\mathrm{d}H}{\mathrm{d}t}=\frac{\partial H}{\partial t}$, hence
\begin{gather*}
  \frac{\mathrm{d}}{\mathrm{d}t}\left(\frac{1}{2}\langle \mathcal{A}(t,q)\dot q,\dot
  q\rangle +\Psi   (t,q)\right)=\frac{\mathrm{d}H(t,q,\dot q)}{\mathrm{d}t}+\frac{\mathrm{d}\Phi(t,q)}{\mathrm{d}t}+
  \frac{\mathrm{d}\Psi (t,q)}{\mathrm{d}t}=\\
 \frac{1}{2}\left\langle \frac{\partial \mathcal{A}(t,q)}{\partial t}\dot q,\dot q\right\rangle+
  \frac{\partial (\Phi (t,q)+\Psi(t,q))}{\partial q_{i}}\dot q_i+\frac{\partial \Psi (t,q)}{\partial t}.
\end{gather*}
By assumption~$(\gamma)$, if $V(t,q,\dot q)\ge 0$, then
\begin{gather*}
  \left|\dot V_{f}(t,q,\dot q)\right|=\\\left(\frac{1}{2}\langle \mathcal{A}(t,q)\dot q,\dot
  q\rangle +\Psi   (t,q)\right)^{-\theta}\left|\frac{\mathrm{d}}{\mathrm{d}t}\left(\frac{1}{2}\langle \mathcal{A}(t,q)\dot q,\dot
  q\rangle +\Psi   (t,q)\right)\right|\le\\
  K\left(\frac{1}{2}\langle \mathcal{A}(t,q)\dot q,\dot
  q\rangle +\Psi (t,q)\right)\le \frac{K}{\kappa}\dot W_f(t,q,\dot q).
\end{gather*}
\end{proof}

\begin{lemma}\label{lem:boundsW} For the  functions $V(\cdot)$ and $W(\cdot)$ defined
defined by \eqref{eq:defofVW}, the  corresponding functions $w_0(\cdot)$, $w^0(\cdot)$
defined by \eqref{eq:w0w0} satisfy the following estimates:
\begin{gather}\label{eq:neww0w01}
w_0(t)\ge \tilde w_0(t):=\min_{q\in \Psi_t^{-1}\left([0,R]\right)}\left\{\langle
a(t,q),q\rangle-\sqrt{2[R-\Psi (t,q)]\langle \mathcal{A}(t,q)q,q\rangle}\right\},\\ \label{eq:neww0w02}
w^0(t)\le \tilde w^0(t):=\max_{q\in \Psi_t^{-1}\left([0,R]\right)}\left\{\langle a(t,q),q\rangle+
\sqrt{2[R-\Psi(t,q)]\langle \mathcal{A}(t,q)q,q\rangle}\right\},\\\label{eq:neww0w03}
 \omega_0:=\inf_{t\in \mathbb R}w_0(t)\ge -\max_{s\in [0,R]}\sqrt{\overline{\Theta}(s)}\left[\sqrt{2(R-s)}+\Xi(s)\right],\\
 \omega^0:=\sup_{t\in \mathbb R}w^0(t)\le \max_{s\in [0,R]}\sqrt{\overline{\Theta}(s)}\left[\sqrt{2(R-s)}+\Xi(s)\right].\label{eq:neww0w04}
\end{gather}
 \end{lemma}
\begin{proof}We know that
\begin{gather*}
  V_t^{-1}(0)=\left\{(q,\dot q)\in \mathbb R^{2m}:\langle \mathcal{A}(t,q)\dot q,\dot q\rangle =2\left[R-\Psi (t,q)\right],
  \;\Psi_t(q)\le R\right\}
\end{gather*}
and  $$\langle a(t,q),q\rangle-\left|\langle \mathcal{A}(t,q)\dot q,q\rangle \right|\le W(t,q,\dot q)\le \langle
a(t,q),q\rangle+\left|\langle \mathcal{A}(t,q)\dot q,q\rangle \right|.$$ By assumptions~$(\beta)$ and $(\delta)$ we
have
\begin{gather*}
  \left|\langle a(t,q),q\rangle\right| \le \Xi(\Psi (t,q))\sqrt{\overline{\Theta}(\Psi (t,q))}.
\end{gather*}
Now to obtain the required estimates it is sufficiently to observe that
\begin{gather*}
\left|\langle \mathcal{A}(t,q)\dot q,q\rangle \right|\left|_{V_t^{-1}(0)}\right.\le \sqrt{\langle \mathcal{A}(t,q)\dot q,\dot
q\rangle\langle \mathcal{A}(t,q)q,q\rangle}\left|_{V_t^{-1}(0)}\right.\le \\ \sqrt{2[R-\Psi
(t,q)]\overline{\Theta}(\Psi (t,q))}
\end{gather*}
and $\Psi (t,q)\ge 0$.
\end{proof}

Now we are in position to prove the following theorem.
\begin{theorem}\label{th:bsLs}
Let for the Lagrangian \eqref{eq:Ltqq} the assumptions~($\alpha$)--($\delta $) be
valid. Then the corresponding Lagrangian system has  a global solution $q_*(t)$
which for some positive number $\sigma \in (0,(\omega^0-\omega_0)/(\kappa R)]$  satisfies the
inequalities
\begin{gather*}
  \tfrac{1}{2}\left\langle \mathcal{A}\left(t,q_*(t)\right)\dot q_*(t),\dot q_*(t)\right\rangle
  +\Psi \left(t,q_*(t)\right)\le \\ \mathfrak f_{\theta,R} \left(\tfrac{K}{2\kappa}\left[\sup_{t\le s\le t+\sigma}\tilde w^0(s)-
  \inf_{t-\sigma \le s\le t}\tilde w_0(s)\right]\right),\\
  \omega_0\le \left\langle \mathcal{A}(t,q)\dot q+a(t,q),q\right\rangle\le \omega^0
\end{gather*}
where
\begin{gather*}
  \mathfrak f_{\theta,R}(z):=  \begin{cases}
    R\,e^z & \text{if}\quad \theta =1, \\
    \left[(1-\theta)z+R^{1-\theta}\right]^{\frac{1}{1-\theta}} & \text{if}\quad \theta \in
    [0,1),
  \end{cases}
\end{gather*}
and the functions $\tilde w_0(t)$, $\tilde w^0(t)$ and numbers $\omega_0,\;\omega^0$ are
defined by \eqref{eq:neww0w01}--\eqref{eq:neww0w04}.
\end{theorem}

\begin{proof}

Let $w_*<\omega_0$ and $w^*> \omega^0$ be arbitrary numbers, where $\omega_0,\omega^0$ are
defined by \eqref{eq:neww0w03},\eqref{eq:neww0w04}. The function $W(\cdot)$ in new
coordinates $q,\;p:=\mathcal{A}(t,q)\dot q+a(t,q)$ takes the form of an indefinite
nondegenerate quadratic form $\langle p,q\rangle $. From this it follows that the set
$\mathcal{W}:=W^{-1}\left((w_*,w^*)\right)$ is connected and for each $t\in \mathbb R$ the function
$W_t(\cdot)$ restricted to the set $V_{t}^{-1}\left((-\infty,0]\right)$ takes its maximal and
minimal values on the boundary $V_t^{-1}(0)$. Hence,
$V^{-1}\left((-\infty,0]\right)\subset \mathcal{W} $ and by Lemma~\ref{lem:VWpair} the
conditions~(A) and (B) are valid with $c_*=c^*=K/\kappa $ and  $\alpha(t)\ge \kappa R$
respectively. Obviously, the functions $\tau_{\pm}(t)$ defined in
Theorem~\ref{th:exbs} satisfy in our case the inequalities
\begin{gather*}
  \left|\tau_{\pm}(t)-t\right|\le \frac{\omega^0-\omega_0}{\kappa R}.
\end{gather*}

Now we define the set
\begin{gather*}
  \mathcal{M}_t:=\{(q,\dot q)\in \mathbb R^{2m}:\dot q=q-\mathcal{A}^{-1}(t,q)a(t,q),\;\langle \mathcal{A}(t,q)q,q\rangle \le w^*\}.
\end{gather*}
Obviously, $0\le W_t(q,\dot q)\left|_{\mathcal{M}_t}\right.=\langle \mathcal{A}(t,q)q,q\rangle \le w^*$.
Since by assumption~($\beta $) $\langle \mathcal{A}(t,q)q,q\rangle $ is spatially coercive, and
\eqref{eq:quasiconvA} implies that
\begin{gather}\label{eq:quasiconvA1}
  \frac{\partial \langle \mathcal{A}(t,q)q,q\rangle}{\partial q_i}q_i=
  \left\langle \left(2\mathcal{A}(t,q)+\frac{\partial \mathcal{A}(t,q)}{\partial q_i}q_i\right)q,q\right\rangle \ge
  \frac{\kappa}{2} \langle \mathcal{A}(t,q)q,q\rangle,
\end{gather}
then $\langle \mathcal{A}(t,q)q,q\rangle $ is regular spatially coercive. For this reason,
$\mathcal{M}_t$ is a compact manifold with boundary.

In order to show that $\nu $ defined by \eqref{eq:defnu} is bounded, note that
in view of assumption~($\beta $) the function $\Psi (t,q)$ is bounded from above by
the constant $\underline{\Theta}^{-1}(w^*)$ on the set where $\langle \mathcal{A}(t,q)q,q\rangle
\le w^*$, and now, taking into account the definition of $V(\cdot)$, it is
sufficient to prove that
\begin{equation}\label{eq:boundnu}
  \sup_{t\in \mathbb R}\max\left\{\langle \mathcal{A}(t,q)\dot q,\dot q\rangle :\dot
q=q-\mathcal{A}^{-1}(t,q)a(t,q),\;\langle \mathcal{A}(t,q)q,q\rangle \le w^*\right\}<\infty.
\end{equation}
But  from ($\delta $) for such points that $\langle \mathcal{A}(t,q)q,q\rangle \le w^*$, we obtain

\begin{gather*}|\langle a(t,q),q\rangle|\le \sqrt{w^*}\Xi(\underline{\Theta}^{-1}(w^*)), \\ \sqrt{\langle
\mathcal{A}^{-1}(t,q)a(t,q),a(t,q)\rangle}\le \Xi(\underline{\Theta}^{-1}(w^*)).\end{gather*}
Hence,
\begin{gather*}
\left\langle \mathcal{A}(t,q)[q-\mathcal{A}^{-1}(t,q)a(t,q)],[q-\mathcal{A}^{-1}(t,q)a(t,q)]\right\rangle =\\ \langle
\mathcal{A}(t,q)q,q\rangle - 2\langle a(t,q),q\rangle +\langle \mathcal{A}^{-1}a(t,q),a(t,q)\rangle \le \\
w^*+2\sqrt{w^*}\Xi(\underline{\Theta}^{-1}(w^*))+\Xi^2(\underline{\Theta}^{-1}(w^*)),
\end{gather*}
and \eqref{eq:boundnu} is proved.

Let us show that the condition~(C) is valid. Since  in $(q,p)$-coordinates the
function $W(q,p)=\langle p,q\rangle $ does not depend on $t$,  it remains only to prove
that for any fixed $t\in \mathbb R$ the set $$\partial \mathcal{M}_t=\{(q,p)\in
\mathbb R^{2m}:p=\mathcal{A}(t,q)q,\;\langle \mathcal{A}(t,q)q,q\rangle = w^*\}$$ is a retract of
$W_t^{-1}(w^*)=\{(q,p)\in \mathbb R^{2m}:\langle p,q\rangle =w^*\}$. Observe that for any
$q\ne0$ from \eqref{eq:quasiconvA1} we get
\begin{gather*}
  \frac{\mathrm{d}}{\mathrm{d}\tau}e^{2\tau}\left\langle \mathcal{A}\left(t,e^\tau q\right)q,q\right\rangle \ge
  \tfrac{\kappa}{2} e^{2\tau}\left\langle \mathcal{A}\left(t,e^\tau q\right)q,q\right\rangle.
\end{gather*}
This implies that for any fixed $q$ the mapping $\tau  \mapsto e^{2\tau}\left\langle
\mathcal{A}\left(t,e^\tau q\right)q,q\right\rangle$ is a  diffeomorphism of $\mathbb R$ onto $(0,\infty)$.
Hence, for any $(q,z)\in \mathbb R^{m}\times(0,\infty)$ there exists a unique $\tau(q,z)$ such
that
\begin{gather*}
  e^{2\tau}\left\langle \mathcal{A}\left(t,e^\tau q\right)q,q\right\rangle\left|_{\tau =\tau(q,z)}\right.=z,\quad
  \tau \left( q,\langle \mathcal{A}(t,q)q,q\rangle \right)=0.
\end{gather*}
By the inverse function theorem the mapping $\tau(\cdot):\left(\mathbb R^{m}\setminus
\{0\}\right)\times(0,\infty) \mapsto \mathbb R$ which we have constructed is smooth. Now the
required retraction is defined by the mapping
\begin{gather*}
  q \mapsto e^{\tau(q,\langle p,q\rangle)}q,\quad p
  \mapsto e^{\tau(q,\langle p,q\rangle)}\mathcal{A}\left(t,e^{\tau(q,\langle p,q\rangle)}q\right)q.
\end{gather*}
Now the existence of searched solution $q_*(t)$ follows from
Theorem~\ref{th:exbs}.
\end{proof}

\begin{corollary}\label{cor:boundedLagr}If the assumptions $(\alpha)$--$(\delta)$ are valid with $\Psi(\cdot)=\Phi(\cdot)$,
then the solution $q_*(t)$ has the property  $\sup_{t\in \mathbb R}\left|L(t,q_*(t),\dot
q_*(t))\right|<\infty $.
\end{corollary}

The next two lemmas will be useful for verifying the assumptions~$(\alpha)$ and
$(\gamma)$.

\begin{lemma}\label{lem:quasiconv} Let there exist positive constants
$\kappa,R_0,c_1,c_2$ such that
\begin{gather}\label{eq:condAtq}
  \min_{\|y\|=1}\frac{\left\langle \left(\mathcal{A}(t,q)+
  \frac{1}{2}\frac{\partial \mathcal{A}(t,q)}{\partial
  q_i}q_i\right)y,y\right\rangle}{\langle \mathcal{A}(t,q)y,y\rangle}\ge \kappa>0\quad \forall(t,q)\in
  \mathbb R^{1+m},
\end{gather}
\begin{gather}\label{eq:dFdqiqi}
 \dfrac{\partial \Phi(t,q)}{\partial q_i}q_i\ge \kappa \Psi (t,q)+\frac{1}{2\kappa}\max_{\|y\|=1}
  \frac{\left\langle a(t,q)+\frac{\partial a(t,q)}{\partial q_i}q_i,y\right\rangle^2}{\langle \mathcal{A}(t,q)y,y\rangle}
  \\ \nonumber
  \text{once}\quad \Psi  (t,q)\ge R_0,
\end{gather}
\begin{gather}\label{eq:dadqiqi}
  \dfrac{\partial \Phi(t,q)}{\partial q_i}q_i\ge -c_1,\quad
  \max_{\|y\|=1}\frac{\left|\left\langle\frac{\partial a(t,q)}{\partial q_i}q_i,y\right\rangle\right|}{\sqrt{\langle
\mathcal{A}(t,q)y,y\rangle}}\le c_2\\ \nonumber
  \text{once}\quad \Psi  (t,q)\le R_0;
\end{gather}

Then  under the assumption~($\delta $), the assumption~$(\alpha)$ is valid with
$$R:=R_0+\left[\frac{\sqrt{2}(c_2+\Xi (R_0))+\sqrt{2(c_2+\Xi (R_0))^2+4\kappa (c_1+\kappa
R_0)}}{2\kappa}\right]^2.$$
\end{lemma}
\begin{proof}From \eqref{eq:condAtq} it follows that
\begin{gather*}
  \dfrac{\partial L}{\partial q_i}q_i+\dfrac{\partial L}{\partial \dot q_i}\dot q_i\ge \kappa \langle \mathcal{A}(t,q)\dot q,\dot q\rangle +
  \left\langle a(t,q)+\dfrac{\partial a(t,q)}{\partial q_i}q_i,\dot q\right\rangle+\dfrac{\partial
\Phi(t,q)}{\partial q_i}q_i
\end{gather*}
If we put $y=\|\dot q\|^{-1}\dot q\in \mathbb S_1(0):=\{y\in \mathbb R^m:\|y\|=1\}$ and
$z:=\sqrt{\langle \mathcal{A}(t,q)\dot q,\dot q\rangle}/\sqrt{2}$, then it is sufficient to show
that the inequality
\begin{gather*}
  z^2+ \Psi  (t,q)\ge R
\end{gather*}
yields
\begin{gather*}
  \kappa z^2-\sqrt{2}
  \frac{\left|\left\langle a(t,q)+\frac{\partial a(t,q)}{\partial q_i}q_i,y\right\rangle\right|}{\sqrt{\langle \mathcal{A}(t,q)y,y}\rangle}z+
  \dfrac{\partial \Phi(t,q)}{\partial q_i}q_i - \kappa  \Psi (t,q)\ge 0.
\end{gather*}
But if $\Psi(t,q)\ge R_0$ then in view of \eqref{eq:dFdqiqi} the  quadratic
polynomial (with respect to $z$) in the left-hand side of the last inequality
takes only nonnegative values for all $y\in \mathbb S_1(0)$. And if $\Psi(t,q)\le R_0$,
then taking into account assumptions \eqref{eq:dadqiqi}, ($\delta$), it is no hard
to show that the greatest root of this polynomial (if it exists) does not
exceed $\sqrt{R-R_{0}}\le \sqrt{R-\Psi(t,q)}$ for all $y\in \mathbb S_1(0)$. Hence, in
this case, the polynomial also takes nonnegative values for $z\ge
\sqrt{R-\Psi(t,q)}$.
\end{proof}

\begin{lemma}\label{lem:guidW}
Let there exist a number $\theta \in [0,1]$ and nonnegative numbers $c_3,\ldots c_8$
such that
 \begin{gather*}
 \max_{\|y\|=1} \frac{\left|\left\langle \frac{\partial \mathcal{A}(t,q)}{\partial t}y,y\right\rangle \right|}{\langle \mathcal{A}(t,q)y,y\rangle}\le
  c_3\Psi^\theta(t,q) +c_4\quad \forall(t,q)\in \mathbb R^{1+m}
\end{gather*}
\begin{gather*}
  \max_{\|y\|=1}\frac{\left|\frac{\partial (\Phi (t,q)+\Psi(t,q))}{\partial q_i}y_i\right|}{\sqrt{\langle \mathcal{A}(t,q)y,y\rangle}}\le
  c_5\Psi^{\theta +1/2}(t,q)+c_6\quad  \forall(t,q)\in \mathbb R^{1+m}
\end{gather*}
\begin{gather*}
 \left|\dfrac{\partial \Psi (t,q)}{\partial t}\right|\le c_7\Psi^{\theta +1}(t,q)+c_8\quad  \forall(t,q)\in
 \mathbb R^{1+m}.
\end{gather*}
Then the assumption~$(\gamma)$ is valid with
\begin{gather*}
  K= c_3+\sqrt{2}c_5+c_7+R^{-\theta}\left( c_4+\sqrt{2}R^{-1/2}c_6+R^{-1}c_8\right).
\end{gather*}

\end{lemma}
\begin{proof} Let again $z:=\sqrt{\langle \mathcal{A}(t,q)\dot q,\dot q\rangle}/\sqrt{2}$. Then we
have
\begin{gather*}
  \left|  \frac{1}{2}\left\langle \frac{\partial \mathcal{A}(t,q)}{\partial t}\dot q,\dot q\right\rangle+
  \frac{\partial (\Phi (t,q)+\Psi(t,q))}{\partial q_{i}}\dot q_i+\frac{\partial \Psi (t,q)}{\partial t}\right|\le\\
  \left(c_3\Psi^\theta (t,q)+c_4\right)z^2+\sqrt{2}\left(c_5\Psi^{\theta +1/2}(t,q)+c_6\right)z+
  c_7\Psi^{\theta
  +1}(t,q)+c_8\le \\
\left(c_3(z^2+\Psi(t,q))^\theta +c_4\right)(z^2+\Psi(t,q))+\\
\sqrt{2}\left(c_5(z^2+\Psi(t,q))^{\theta +1/2}+c_6\right)\sqrt{z^2+\Psi(t,q)}+\\
c_7(z^2+\Psi(t,q))^{\theta
  +1}+c_8\le K(z^2+\Psi(t,q))^{\theta +1}.
\end{gather*}
\end{proof}
Let us now discuss the uniqueness problem. Usually, to guarantee the uniqueness
of bounded solutions (in particular, almost periodic solutions) to Lagrangian
systems, the convexity of Lagrang\-i\-an function is required. In Cieutat's
paper \cite{Cie03b} it is assumed that the function $\frac{\partial L(t,\cdot)}{\partial
u}:\mathbb R^{2m} \mapsto \mathbb R^{2m}$ is globally Lipschitzian with time independent
Lipschitz constant, and the convexity condition is formulated as follows: there
exists a constant $c>0$ such that
\begin{gather}\label{eq:str_conv}
  \left(\dfrac{\partial L(t,u)}{\partial u_i}-\dfrac{\partial L(t,v)}{\partial
  v_i}\right)(u_i-v_i)\ge c\|u-v\|^2\quad \\
  \nonumber \forall u:=(q',\dot q'),\;v:=(q'',\dot q'')\in
  \mathbb R^{2m}.
  \end{gather}
It should be noted that for Lagrangian~\eqref{eq:Ltqq}, in the case where
$\mathcal{A}(t,q)$ nonlinearly depends on $q$, the above global conditions look
unnatural (see the Remark~\ref{rem:restrF} below).

For Lagrangian~\eqref{eq:Ltqq}, we are going to relax the conditions of
\cite{Cie03b}  via the Theorem~\ref{th:unbs}. (However, here for simplicity we
consider the case where $\mathcal{A}(\cdot),\;a(\cdot)$ and $\Phi(\cdot)$ are
$\mathrm{C}^{2}$-mappings).

Put
\begin{gather*}
  \tilde V(t,u,v):=\max\left\{\frac{1}{2}\langle \mathcal{A}(t,q')\dot q',\dot q'\rangle
+\Psi (t,q'),\frac{1}{2}\langle \mathcal{A}(t,q'')\dot q'',\dot q''\rangle +\Psi (t,q'')\right\},
\end{gather*}
and denote
\begin{gather*}
  N(t;r):=\sup \left\{\frac{\left\|\tfrac{\partial L(t,u)}{\partial \dot q'}-\tfrac{\partial L(t,v)}{\partial \dot
  q''}\right\|}{\|u-v\|}: (u,v)\in \tilde V_t^{-1}\left((-\infty,r]\right), u\ne v\right\}.
\end{gather*}
Let $\lambda(t,q)$ and $\Lambda(t,q)$ be, respectively, the minimal and the maximal
eigenvalues of operator $\mathcal{A}(t,q)$. Define
\begin{gather*}
  \vartheta(t;r):= \max\left\{\frac{\Lambda(t,q'')}{\lambda(t,q')}:q'\in \Psi_t^{-1}([0,r]),
  q''\in \Psi_t^{-1}([0,r])\right\}.
\end{gather*}

For any set $\tilde \Omega \subset \mathbb R^{1+2m}$ we  define the set
\begin{gather*}
    \tilde \Omega^*:=\left\{(t,u,v)\in \mathbb R^{1+4m}:(t,u)\in \tilde \Omega,\;(t,v)\in \tilde \Omega \right\}
  \end{gather*}
(see Theorem~\ref{th:unbs}).

\begin{theorem}\label{th:uniqLagr}Let the assumptions $(\beta)$ and $(\delta)$  be
valid and let for a set $\tilde \Omega \subset \mathbb R^{1+2m}$ there exist  numbers $r>0$
and $d\ge 1$ such that $\tilde \Omega^* \subseteq \tilde V^{-1}\left([0,r]\right)$ and
\begin{gather*}
 \varrho(t;r,d):= \inf\left\{\frac{\left(\tfrac{\partial L(t,u)}{\partial u_i}-\tfrac{\partial L(t,v)}{\partial
  v_i}\right)(u_i-v_i)}{\|u-v\|^{2d}}: (u,v)\in \tilde \Omega^*_t,\;u\ne v\right\}>0.
\end{gather*}
 Suppose in addition that
  \begin{gather*}
  I(t;r,d):=\left|\int_{0}^{t}\frac{\varrho(s;r,d)}{N^d(s;r)}\,\mathrm{d}s\right|\to \infty,\quad t\to \pm \infty
\end{gather*}
and if $d=1$, then also
\begin{gather*}
  \liminf_{t\to \pm \infty}\frac{\ln \left(1+\sqrt{\vartheta(t;r)}\right)}{I(t;r,1)}<1.
\end{gather*}
Then the Lagrangian system cannot have  two different global solutions
$q_j(t),t\in \mathbb R$, j=1,2, such that $(t,q_j(t),\dot q_j(t))\in \tilde \Omega $ for all
$t\in \mathbb R$, $j=1,2$.
\end{theorem}

\begin{proof}
In order to apply the Theorem~\ref{th:unbs}, we introduce the function
\begin{gather*}
  U(t,u,v):=\left(\dfrac{\partial L(t,u)}{\partial \dot q'_i}-\dfrac{\partial L(t,v)}{\partial \dot
  q''_i}\right)(q_i'-q_i'')=\\
  \langle \mathcal{A}(t,q')\dot q'+a(t,q')-\mathcal{A}(t,q'')\dot q''-a(t,q''),q'-q''\rangle.
\end{gather*}

In the same way as in Lemma~\ref{lem:boundsW}, one can show that
\begin{gather*}
|\langle \mathcal{A}(t,q)\dot q+a(t,q),q\rangle|\le \sqrt{\overline{\Theta}(\Psi (t,q))}\left[\sqrt{2[r-\Psi
(t,q)]}+\Xi\left(\Psi (t,q)\right)\right]\\ \forall(t,q,\dot q)\in \tilde \Omega,
\end{gather*}
and since $\langle \mathcal{A}(t,q')q'', q''\rangle\le \frac{\Lambda(t,q')}{\lambda
(t,q'')}\overline{\Theta}\left(\Psi (t,q'')\right)$, then
\begin{gather*}
|\langle \mathcal{A}(t,q')\dot q'+a(t,q'),q''\rangle|\le\\
  \sqrt{\langle \mathcal{A}(t,q')q'',q''\rangle}\left[\sqrt{\langle \mathcal{A}(t,q')\dot q',\dot q'\rangle}+\Xi\left(\Psi (t,q')\right)\right]
  \le\\
  \sqrt{ \vartheta(t;r)\overline{\Theta}(\Psi (t,q''))}\left[\sqrt{2[r-\Psi (t,q')]}+\Xi\left(\Psi (t,q')\right)\right]\quad \forall(t,u,v)\in \tilde \Omega^*.
\end{gather*}

 Now we have
\begin{gather}\label{eq:estUtuv}
|U(t,u,v)|\le 2\omega^*(r)\left[1+\sqrt{\vartheta(t;r)}\right]\quad  \forall(t,u,v)\in \tilde \Omega^*
\end{gather}
where
\begin{gather*}
  \omega^*(r):=\sqrt{\overline{\Theta}(r)}\max_{1\le s\le r}\left[\sqrt{2[r-s]}+\Xi(s)\right].
\end{gather*}

Hence,  in the case under consideration, the function $b(t,r)$ from
Theorem~\ref{th:unbs} satisfies the inequality
\begin{gather*}
  b(t,r)\le 2\omega^*(r)\left(1+\sqrt{\vartheta(t;r)}\right).
\end{gather*}

Nextly, the inequality
\begin{gather*}
|U(t,u,v)|\le N(t;r)\|u-v\|^2
\end{gather*}
together with conditions imposed on $L(\cdot)$ yields
\begin{gather*}
  \dot U(t,u,v)_{(f,f)}=
  \left(\dfrac{\partial L(t,u)}{\partial u_i}-\dfrac{\partial L(t,v)}{\partial
  v_i}\right)(u_i-v_i)\ge \varrho(t;r,d)\|u-v\|^{2d}\ge\\ \beta(t,r)|U(t,u,v)|^{d}
\end{gather*}
if $(t,u,v)\in \tilde \Omega^*$ where $\beta(t;r,d):=\varrho(t;r,d)N^{-d}(t,r)$. Now if
we put $h(u)=\int_{1}^{u}s^{-d}\,\mathrm{d}s$, then the reasoning which we used when
proving the Theorem~\ref{th:unbs} yields the assertion of the
Theorem~\ref{th:uniqLagr}.
 \end{proof}

 It appears that  instead of the convexity condition of Theorem~\ref{th:uniqLagr} it is preferable to verify
an analogous assumption for corresponding Hamiltonian
\begin{gather}\label{eq:corrHamfunc}
  H(t,z)\equiv H(t,q,p):=
\frac{1}{2}\left\langle \mathcal{A}^{-1}(t,q)(p-a(t,q)),p-a(t,q)\right\rangle -\Phi(t,q)\\ (z:=(q,p)).
\nonumber \end{gather}

Let us introduce the function
\begin{gather*}
  Y(t,z):=\frac{1}{2}\left\langle\mathcal{A}^{-1}(t,q)(p-a(t,q)),p-a(t,q)\right\rangle+\Psi(t,q),
\end{gather*}
which corresponds to the function $\tfrac{1}{2}\left\langle\mathcal{A}(t,q)\dot q,\dot q
\right\rangle+\Psi(t,q)$.

Let $\mathrm{Id}_m$ and $0_m$  be the identity matrix and the zero matrix of
dimensions $m$ respectively. Introduce the matrices
\begin{gather*}
  I:=\begin{pmatrix}
    -\mathrm{Id}_m & 0_m \\
    0_m & \mathrm{Id}_m
  \end{pmatrix}, \quad J:=
  \begin{pmatrix}
    0_m & \mathrm{Id}_m\\
        -\mathrm{Id}_m &  0_m
  \end{pmatrix}
\end{gather*}
Put $z':=(q',p'),\;z'':=(q'',p'')$ and denote by $\hat  V(t,z',z'')$ the function
obtained from $\tilde  V(t,u,v)$ after the substitutions $\dot q'=
\mathcal{A}^{-1}(t,q')(p'-a(t,q'))$, $\dot q''= \mathcal{A}^{-1}(t,q'')(p''-a(t,q''))$.
Obviously,
\begin{gather*}
  \hat V(t,z',z''):=\max\left\{Y(t,z'),Y(t,z'')\right\}
\end{gather*}

\begin{theorem}\label{th:uniqHam}
Let the assumptions $(\beta)$ and $(\delta)$  be valid and let for a set $\hat \Omega
\subset \mathbb R^{1+2m}$ there exist  numbers $r>0$ and $d\ge 1$ such that $\hat \Omega^*
\subseteq \hat  V^{-1}\left([0,r]\right)$ and
\begin{gather*}
 \hat  \varrho(t;r,d):= \inf\left\{
\tfrac{\left\langle I\left(\tfrac{\partial H(t,z')}{\partial z'}-\tfrac{\partial H(t,z'')}{\partial
z''}\right),z'-z''\right\rangle}{\|z'-z''\|^{2d}}:
  (z',z'')\in \hat \Omega^*_t,\;z'\ne z''\right\}>0.
\end{gather*}
Suppose in addition that
 $\lim_{t\to \pm \infty} \left|\int_{0}^{t}\hat \varrho(s;r,d)\,\mathrm{d}s\right|=\infty $ and if $d=1$, then also
\begin{gather*}
  \liminf_{t\to \pm \infty}\frac{\ln \left(1+\sqrt{\vartheta(t;r)}\right)}{2\left|\int_{0}^{t} \hat
\varrho(s;r,d)\,\mathrm{d}s\right|}<1.
\end{gather*}
Then the system with Hamiltonian \eqref{eq:corrHamfunc} cannot have two
different global solutions $\left(q_j(t),p_j(t)\right),t\in \mathbb R$, j=1,2, such that
$\left(t,q_j(t),p_j(t)\right)\in \hat \Omega $ for all $t\in \mathbb R$, $j=1,2$.
\end{theorem}

\begin{proof}
In order to apply the Theorem~\ref{th:unbs} in the case of Hamiltonian system
\begin{gather*}
  \dot z=JH'_z(t,z),
\end{gather*}
introduce the function $\hat U(z',z''):=\langle q'-q'',p'-p''\rangle $. After the
substitutions $p'=\frac{\partial L(t,u)}{\partial \dot q'}$, $p''=\frac{\partial L(t,v)}{\partial \dot
q''}$, this function coincides with the function $U(t,u,v)$ which appears when
proving the Theorem~\ref{th:uniqLagr}. Hence, the estimate \eqref{eq:estUtuv}
implies that
\begin{gather*}
  \left|\hat U(z',z'')\right|\le 2\omega^*(r)\left[1+\sqrt{\vartheta(t;r)}\right]
\end{gather*}
once $\hat V(t,z',z'')\le r$, and the inequality
\begin{gather*}
|\hat U(z',z'')|\le \frac{1}{2}\|z'-z''\|^2
\end{gather*}
together with definition of $\hat \varrho(t;r,d)$ yields
\begin{gather*}
  \dot {\hat U}(z',z'')_{(JH'_{z},JH'_{z})}=
 \left\langle I\left(\dfrac{\partial H(t,z')}{\partial z'}-\dfrac{\partial H(t,z'')}{\partial z''}\right),z'-z''\right\rangle
  \ge \\ \hat \varrho(t;r,d)\|z'-z''\|^{2d}\ge  2^{d} \hat \varrho(t;r,d)|\hat   U(z',z'')|^d
\end{gather*}
if $\hat  V(t,z',z'')\le r$. The rest of the proof is based on the same arguments
as the proof of previous theorem.
\end{proof}

As a corollary of Theorems~\ref{th:bsLs},~\ref{th:uniqHam} we can  get new
sufficient conditions for the existence of almost periodic solutions to
Lagrangian systems. Namely, consider the case where the following assumption is
valid:

\begin{itemize}
  \item[$(\epsilon)$] the mappings $\mathcal{A}(\cdot,q):\mathbb R \mapsto \mathrm{Hom}(\mathbb R^{m})$, $a(\cdot,q):\mathbb R \mapsto \mathbb R^{m}$,
  $\Phi(\cdot,q):\mathbb R \mapsto \mathbb R$ together with their first order partial derivatives in
  $q$ are almost periodic uniformly for $q\in \mathbb R^m$ and the function
  $\Psi_*(q):=\inf_{t\in \mathbb R}\Psi(t,q)$ is coercive.
\end{itemize}

Denote
\begin{gather*}
  \Lambda^*(q):=\sup_{t\in \mathbb R}\Lambda(t,q),\quad \alpha_*(q):=\sup_{t\in \mathbb R}\|a(t,q)\|.
\end{gather*}
Since
\begin{gather*}
  Y(t,z)\ge \frac{1}{2\Lambda^*(q)}\left(\|p\|-\alpha_*(q)\right)^2+\Psi_*(q),
\end{gather*}
and the function in the right-hand side of this inequality is coercive, then
for any $r>0$ the set
\begin{gather*}
  \mathcal{V}(r):=\mathrm{cls}\bigcup_{t\in \mathbb R} \left\{(p,q)\in \mathbb R^{2m}: Y(t,p,q)\le r,\;
  \omega_0\le\langle p,q\rangle \le\omega^0\right\}
\end{gather*}
is compact (see \eqref{eq:neww0w03},\eqref{eq:neww0w04} for definitions of
$\omega_0$, $\omega^0$).

\begin{theorem}\label{th:apsolLagr} Let the assumptions $(\alpha)-(\epsilon)$ be
valid. Put $$r:=\mathfrak f_{\theta,R} \left(\tfrac{K}{2\kappa}(\omega^0-\omega_0)\right)$$ (the function
$\mathfrak f_{\theta,R}(\cdot)$ is defined in Theorem~\ref{th:bsLs}) and suppose that there
exist numbers $\varrho_*>0$ and $d\ge 1$ such that
\begin{gather*}
  \left\langle I\left(\dfrac{\partial H(t,z')}{\partial z'}-\dfrac{\partial H(t,z'')}{\partial
  z''}\right),z'-z''\right\rangle\ge \varrho_*\|z'-z''\|^{2d}
\end{gather*}
for all $(t,z',z'')\in \mathbb R\times \mathcal{V}(r)\times \mathcal{V}(r)$. Then the set $\mathcal{V}(r)$
contains one and only one global solution of the system with Hamiltonian
$H(t,z)$, and this solution is almost periodic.
\end{theorem}

\begin{proof}  By Theorem~\ref{th:bsLs} for any $s\in \mathbb R$,
the Hamiltonian system
\begin{equation}\label{eq:HamSyss}
  \dot z=JH'_z(t+s,z)
\end{equation}
has a global solution taking values in $\mathcal{V}(r)$. Moreover, the same reasoning
as in the proof of Theorem~\ref{th:uniqHam} shows that the set $\mathcal{V}(r)$
contains no other global solutions of system~\eqref{eq:HamSyss}. Now to
complete the proof, it remains only to apply the Amerio theorem (see, e.g.,
\cite{Fin74}).
\end{proof}

Observe now that under the conditions imposed on $L(\cdot)$ the Hamiltonian belongs
to $\mathrm{C}^{2}(\mathbb R^{1+2m}\!\mapsto \!\mathbb R)$. If we denote by $H''_{qq}(t,z)$ the partial
Hesse matrix $\left\{\frac{\partial H(t,q,p)}{\partial q_iq_j}\right\}_{i,j=1}^m$, then it is
easily seen that
\begin{gather*}
\left\langle I\left(\dfrac{\partial H(t,z')}{\partial z'}-\dfrac{\partial H(t,z'')}{\partial
  z''}\right),z'-z''\right\rangle=\\
 \left\langle \left[\int_0^1\mathcal{A}^{-1}(t,sq'+(1-s)q'')\,\mathrm{d}s\right](p'-p''),p'-p''\right\rangle
 -\\
  \left\langle \left[\int_{0}^{1}H''_{qq}(t,sz'+(1-s)z'')\,\mathrm{d}s\right](q'-q''),q'-q''\right\rangle
.
\end{gather*}

Since the first summand of the right-hand side of this equality is positive
definite quadratic form with respect to $p'-p''$ we arrive at conclusion that
for the case where $d=1$, in order that the function $\hat \varrho(t;r,d)$  be
well-defined and positive, it is necessary that
\begin{gather}\nonumber
 \min_{\|\eta \|=1} \left\{-\frac{1}{2}\left[\frac{\partial^2\langle \mathcal{A}^{-1}(t,q)p,p \rangle}{\partial q_i\partial
  q_j}\eta_i\eta_j\right]+\left[\frac{\partial^2\langle \mathcal{A}^{-1}(t,q)a(t,q),p\rangle}{\partial q_i\partial
  q_j}\eta_i\eta_j\right]-\right.\\ \left.\frac{1}{2}\frac{\partial^2\langle \mathcal{A}^{-1}(t,q)a(t,q),a(t,q)\rangle}{\partial q_i\partial
  q_j}\eta_i\eta_j+\frac{\partial^2\Phi(t,q)}{\partial q_i\partial q_j}\eta_i\eta_j\right\}>0 \label{eq:sufconvH}
\end{gather}
for all  $(t,q,p)\in \hat \Omega $, and it is sufficient that the last inequality
holds for all $(t,q,p)$ such that $t\in \mathbb R$ and $(p,q)$ belongs to the convex
hull of the set $\hat \Omega_t $. Observe that  the last set is contained in the
convex hull of the set $Y_t^{-1}([0,r])$. If we treat the left hand side of the
inequality~\eqref{eq:sufconvH} as a quadratic polynomial with respect to
$u=\|p\|$, then we arrive at the following result.

\begin{lemma}\label{lem:convH}
 Put
\begin{gather*}
  \alpha_1(t,q):=1/(2\Lambda(t,q)),\quad \beta_1(t,q):=\|\mathcal{A}^{-1}(t,q)a(t,q)\|,\\
  \gamma_1(t,q):=\langle \mathcal{A}^{-1}(t,q)a(t,q),a(t,q)\rangle +\Psi(t,q),
\end{gather*}
\begin{gather*}
  \alpha_2(t,q):= \max_{\|y\|=1,\|\eta \|=1}\frac{\partial^2\langle \mathcal{A}^{-1}(t,q)y,y \rangle}{\partial q_i\partial
  q_j}\eta_i\eta_j,\\
   \beta_2(t,q):=\max_{\|\eta \|=1}\left\|\frac{\partial^2 \mathcal{A}^{-1}(t,q)a(t,q)}{\partial q_i\partial
  q_j}\eta_i\eta_j\right\|,\\
\gamma_2(t,q):=\min_{\|\eta \|=1}\left[2\frac{\partial^2\Phi(t,q)}{\partial q_i\partial q_j}-\frac{\partial^2\langle
\mathcal{A}^{-1}(t,q)a(t,q),a(t,q)\rangle}{\partial q_i\partial
  q_j}\right]\eta_i\eta_j,
\end{gather*}
and suppose that for any $t\in \mathbb R$ the function $\Psi(t,\cdot):\mathbb R^{m} \mapsto \mathbb R$ is
quasiconvex and that there exists $r>0$ such that for all $(t,q)\in
\Psi^{-1}\left([0,r]\right)$ the inequalities
\begin{gather*}
  \alpha_1(t,q)u^2-2\beta_1(t,q)u+\gamma_1(t,q)\le r,\quad u\ge 0
\end{gather*}
yield the inequality
\begin{gather*}
  \alpha_2(t,q)u^2+2\beta_2(t,q)u-\gamma_2(t,q)<0.
\end{gather*}
Then the inequality \eqref{eq:sufconvH} is valid for all $(t,q,p)$ such that
$t\in \mathbb R$ and $(q,p)$ belongs to convex hull of the set $Y_t^{-1}([0,r])$.
\end{lemma}
\begin{remark}\label{rem:case_a_0}In the particular case where $a(t,q)=0$ the set
$Y_t^{-1}([0,r])$ is convex if for any $t\in \mathbb R$ the function $\Psi(t,\cdot):\mathbb R^{m}
\mapsto \mathbb R$ is quasiconvex.\end{remark}
\begin{remark}\label{rem:restrF}Since for any fixed $t$ and $y\in \mathbb R^m$ the function
$\langle \mathcal{A}^{-1}(t,\cdot)y,y\rangle $ is positive, it cannot be globally strictly concave.
And if $\alpha_2(t,q)>0$ at some point $(t,q)$, then the inequality
\eqref{eq:sufconvH} fails for all $p$ with sufficiently large norm.
\end{remark}

\begin{example}\label{ex:apsLagr}Consider a Lagrangian system which descibes motion
of a particle constrained to move on time-varying helicoid under the impact of
force of gravity and repelling potential field of force.  The vibrating
helicoid is given in 3-D space by the equations
\begin{gather*}
  \mathbf{r}=(q_1\cos q_2, q_1\sin q_2,\chi(t)
  q_2),\quad(q_1,q_2)\in \mathbb R^{2}
\end{gather*}
where $\chi(\cdot)\in \mathrm{C}^{3}(\mathbb R\!\mapsto \!(0,\infty))$ is a given function. Suppose that
the function of repelling potential field is $\Pi(\mathbf{r})=
-k\left(\|\mathbf{r}\|^2+\|\mathbf{r}\|^4\right)$ where $k\ge 1$ is a parameter.
\end{example}
Having assumed for simplicity the mass of particle and the acceleration of
gravity to be unities, we get the following expression for  kinetic energy

$$\frac{1}{2}\lVert\dot{\mathbf{r}}\rVert=\frac{1}{2}\left(\dot{q}_1^2+(\chi^2(t)+q_1^2)\dot{q}_2^2\right)+
\chi(t)\dot{\chi}(t)q_2\dot q_2 +\frac{1}{2} \dot{\chi}^2(t)q_2^2.$$

Since the term $\chi(t)\dot \chi(t)q_2\dot q_2 +\frac{1}{2}\dot \chi^2(t)q_2^2$ gives the
same contribution  into the equations of motion as the term
$-\frac{1}{2}\chi(t)\ddot \chi(t)q_2^2$, we obtain the following Lagrangian
\begin{gather*}
  L(t,q,\dot q)= \frac{1}{2}\left(\dot q_1^2+\left(\chi^2(t)+q_1^2\right)\dot q_2^2\right)-
  \chi(t)q_2-\frac{1}{2}\chi(t)\ddot
  \chi(t)q_2^2+\\ k\left[q_1^2+\chi^2(t)q_2^2+(q_1^2+\chi^2(t)q_2^2)^2\right]
\end{gather*}
 Hence, in this case $a(t,q)=0$,
\begin{gather*}
\langle \mathcal{A}(t,q)\dot q,\dot q\rangle =\dot q_1^2+\left(\chi^2(t)+q_1^2\right)\dot q_2^2,\\
\Phi(t,q)=k\left[q_1^2+\xi(t)q_2^2+(q_1^2+\chi^2(t)q_2^2)^2\right]-\chi(t)q_2
\end{gather*}
where $\xi(t):=\chi^2(t)-\tfrac{1}{2k}\chi(t)\ddot \chi(t)$.

We suppose that the function $\chi(t)$ satisfies the following conditions:
\begin{gather*}
\inf_{t\in \mathbb R}\chi(t)=:\chi_*\ge 1,\quad \sup_{t\in \mathbb R}\left|
\frac{1}{\chi (t)}\frac{\mathrm{d}^i\chi(t)}{\mathrm{d}t^i}\right|=:\eta_i^*<\infty,\;i=1,2,3,\quad \eta_2^*\le
k.
\end{gather*}
Obviously that in this case $\xi(t)>\chi^2(t)/2$.

 Put
\begin{gather*}
  \Psi(t,q)= k \left[q_1^2+\xi(t)q_2^2+2\left(q_1^2+\chi^2(t)q_2^2\right)^2\right].
\end{gather*}
Then
\begin{gather*}
  \dfrac{\partial L}{\partial q_i}q_i+\dfrac{\partial L}{\partial \dot q_i}\dot q_i =
  \langle \mathcal{A}(t,q)\dot q,\dot q\rangle +q_1^2\dot q^2_2+\\ 2k
  \left[q_1^2+\xi(t)q_2^2+2(q_1^2+\chi^2(t)q_2^2)^2\right]-\chi (t)q_2\ge\\
\langle \mathcal{A}(t,q)\dot q,\dot q\rangle+2\Psi(t,q)-\sqrt[4]{\frac{1}{4k}}\sqrt[4]{\langle \mathcal{A}(t,q)\dot
q,\dot q\rangle+2\Psi(t,q)}\ge\\ \kappa(R,k) \left(\tfrac{1}{2}\langle \mathcal{A}(t,q)\dot q,\dot q \rangle +\Psi
(t,q)\right)
\end{gather*}
if $\left(\tfrac{1}{2}\langle \mathcal{A}(t,q)\dot q,\dot q \rangle +\Psi (t,q)\right)\ge R$, where
$\kappa(R,k) :=2-R^{-3/4}(2k)^{-1/4}$. Hence, the assumption $(\alpha)$ holds for
arbitrary $R\ge (32k)^{-1/3}$.

Since $\langle \mathcal{A}(t,q)q,q\rangle =q_1^2+ \left(\chi^2(t)+q_1^2\right)q_2^2$, then the assumption
$(\beta)$ is valid with appropriately chosen function $\underline{\Theta}(\cdot)$ and
with $\overline{\Theta}(\Psi)=2\Psi /k $.

Now let us verify the assumption $(\gamma)$.  We have
\begin{gather*}
  \left|\frac{1}{2}\left\langle \frac{\partial \mathcal{A}(t,q)}{\partial t}\dot q,\dot q\right\rangle\right|\le
  \sup_{t\in \mathbb R}\left|\frac{\dot \chi(t)}{\chi(t)}\right|\langle \mathcal{A}(t,q)\dot q,\dot q\rangle=\eta_1^*\langle \mathcal{A}(t,q)\dot q,\dot q\rangle,
\end{gather*}
\begin{gather*}
  \left|\frac{\partial \Psi (t,q)}{\partial t}\right|\le \sup_{t\in \mathbb R}\max \left\{\left|\frac{\dot
  \xi(t)}{\xi(t)}\right|,4\left|\frac{\dot \chi(t)}{\chi(t)}\right|\right\}\Psi (t,q)\le
  (5\eta_1^*+\eta_3^*)\Psi(t,q),
\end{gather*}
and since $z\le z^3/3+2/3$ for all $z\ge0$ and $\xi(t)/\chi^2(t)\le
1+\eta_2^*/(2k)$, then
\begin{gather*}
 \left| \frac{\partial (\Phi (t,q)+\Psi(t,q))}{\partial q_{i}}\dot q_i\right|\le \\ \sqrt{\langle \mathcal{A}(t,q)\dot q,\dot
 q\rangle}
 \left[4k\sqrt{q_1^2+\left[\frac{\xi(t)}{\chi(t)}\right]^2q_2^2}+12k\left(
 q_1^2+\chi^2(t)q_2^2\right)^{3/2}+1\right]\le \\
 \sqrt{\langle \mathcal{A}(t,q)\dot q,\dot q\rangle}\times \\
  \left(4k\sup_{t\in \mathbb R}\max \left\{1,\frac{\xi (t)}{\chi^2(t)}\right\}
  \sqrt{q_1^2+\chi^2(t)q_2^2}+12k(q_1^2+\chi^2(t)q_2^2)^{3/2}+1\right)\le\\
  \sqrt{\langle \mathcal{A}(t,q)\dot q,\dot q\rangle}\left[ \left(6(2k)^{1/4}+\tfrac{2}{3}\left((2k)^{3/4}+(2k)^{-1/4}\eta_2^*\right)\right)
  \Psi^{3/4}(t,q)+\right. \\ \left.
  \tfrac{4}{3}\left((2k)^{3/4}+(2k)^{-1/4}\eta_2^*\right)+1\right].
\end{gather*}
The same arguments as in the proof of Lemma~\ref{lem:guidW} allows us to assert
that for $\tfrac{1}{2}\langle \mathcal{A}(t,q)\dot q,\dot q\rangle +\Psi (t,q)\ge R$ the assumption
$(\gamma)$ is valid with $\theta =1/4$ and
\begin{gather*}
  K=K(R,k):= \sqrt{2}\left[6(2k)^{1/4}+\tfrac{2}{3}\left((2k)^{3/4}+(2k)^{-1/4}\eta_2^*\right)\right]+ \\
  \left(5\eta_1^*+\eta_3^*\right)R^{-1/4}+
  \sqrt{2}\left[\tfrac{4}{3}\left((2k)^{3/4}+(2k)^{-1/4}\eta_2^*\right)+1\right]
    R^{-3/4}.
\end{gather*}

Lastly, Lemma~\ref{lem:boundsW} yields
\begin{gather*}
  \omega^0\le \max_{s\in [0,R]}\sqrt{4(R-s)s/k}=R\sqrt{2/k},\quad \omega_0\ge -R\sqrt{2/k}.
\end{gather*}
Hence, by Theorem~\ref{th:bsLs}, there exists a global solution $q_*(t),\;t\in
\mathbb R$, satisfying the inequality
\begin{gather*}
  \frac{1}{2}\langle \mathcal{A}(t,q_*(t))\dot q_*(t),\dot q_*(t)\rangle +\Psi(t,q_*(t))\le
  r(k,R)
\end{gather*}
where
\begin{gather*}
r(k,R):=  \left[\frac{3\sqrt{2}R K(R,k)}{4\kappa(R,k)\sqrt{k}}+R^{3/4}\right]^{4/3}.
\end{gather*}
 Observe, that if we put $R=(2k)^{-1/3}$, then $\kappa\left(k,(2k)^{-1/3}\right)=1$ and
\begin{gather*}
  K \left(k,(2k)^{-1/3}\right)\le
  3.78k +1.59k^{3/4}+11.78k^{1/4}+\\ (5.3\eta_1^*+1.1\eta_3^*)k^{1/12}+1,28\eta_2^* \le
  \left(17.15+5.3\eta_1^*+1.28\eta_2^*+1.1\eta_3^*\right)k.
\end{gather*}
Hence,
\begin{gather*}
  r\left(k,(2k)^{-1/3}\right)\le Ck^{2/9}
\end{gather*}
where
\begin{gather*}
  C:= \left(15.28+4.47\eta_1^*+1.08\eta_2^*+0.93\eta_3^*\right)^{4/3}.
\end{gather*}

   From this it
follows that
\begin{gather*}
 \left(q^2_{1*}(t)+\chi^2(t)q^2_{2*}(t)\right)^2\le \tfrac{C}{2}k^{-7/9} \quad \forall t\in \mathbb R,\end{gather*}
 and thus, we obtain the following estimate for the global solution $q_*(t)$:
 \begin{gather*}\|q_*(t)\|^2\le q^2_{1*}(t)+\chi^2(t)q^2_{2*}(t)\le
 \sqrt{\tfrac{C}{2}}k^{-7/18}\quad \forall t\in \mathbb R.
\end{gather*}

Now consider the case where the function $\chi(\cdot)$ is almost periodic together
with its derivatives up to the third order. In order to apply
Lemma~\ref{lem:convH}, observe that
\begin{gather*}
  \alpha_1(t,q)=\frac{1}{2(\chi^2(t)+q_1^2)},\quad
  \alpha_2(t,q)=2\frac{3q^2_1-\chi^2(t)}{(\chi^2(t)+q_1^2)^3},\\
   \beta_1(t,q)=\beta_2(t,q)=0,
\end{gather*}
and it is no hard to show that in our case
\begin{gather*}
  \gamma_2(t,q)\ge 2k \left[\min\{1,\xi(t)\}+2(q_1^2+\chi^2(t)q_2^2)\right],
\end{gather*}
Now it is easily seen that the conditions of Lemma~\ref{lem:convH} will hold
true if on the set where $\Psi(t,q)\le r$ there holds the inequality
\begin{gather*}
 2\frac{3q^2_1-\chi^2(t)}{(\chi^2(t)+q_1^2)^2}(r-\Psi (t,q))\le k
 \left[\min\{1,\xi(t)\}+2(q_1^2+\chi^2(t)q_2^2)\right].
\end{gather*}
Observe that $\sup_{u\ge 0}\frac{3u^2-\chi^2(t)}{(\chi^2(t)+u^2)^2}=\frac{9}{16
\chi^2(t)}\le \frac{9}{16 \chi_*^2}$ and $\xi(t)\ge \chi_*^2/2$. Thus, in the case
where $r=r\left(k,(2k)^{-1/3}\right)$, we get the following sufficient condition for
almost periodicity of solution $q_*(t)$ in terms of restrictions on
parameter~$k$:
\begin{gather*}
  \tfrac{9}{8}Ck^{2/9}\le k\chi_*^2\min\{1,\chi_*^2/2\},\quad k\ge \max\{1,\eta_2^*\}
  \end{gather*}
or
\begin{gather*}
  k\ge \max\left\{1,\eta_2^*,\left[\frac{9C}{8\chi_*^2\min\{1,\chi_*^2/2\}}\right]^{9/7}\right\}.
\end{gather*}

{\bf Conclusions}.

The technique applied in this paper for studying essentially nonlinear
nonauto\-nom\-ous systems  by means of a pair of auxiliary functions allows us
to generalize a number of earlier known results concerning the  questions of
existence and uniqueness of bounded, proper and almost periodic solutions. In
the case where the estimating function is a quadratic form with varying matrix,
the estimates obtained for V-bounded solutions can be efficiently applied to
describe asymptotic behavior of solutions when $t\to \pm \infty $. For Lagrangian
systems  with certain directional quasiconvexity property, there exists a V-W
pair  which allows to establish sufficient conditions for existence of
V-bounded solutions. Our approach yields uniqueness theorems for V-bounded
solutions as well. As a consequence of that, we have obtained new sufficient
conditions for the existence of almost periodic solutions to Lagrangian
systems.

This work was partially supported by the Fundamental Research State Fund of
Ukraine (Project 29.1/025).


\end{document}